\documentclass[a4paper]{amsart}
\usepackage{amsmath,amsthm,amssymb}
\usepackage[dvipdfmx]{graphicx}
\usepackage[all]{xy}
\usepackage{graphicx}
\usepackage{here}
\usepackage{pinlabel}

\newtheorem{theo}{Theorem}[subsection]
\newtheorem{defi}[theo]{Definition}
\newtheorem{lem}[theo]{Lemma}
\newtheorem{rem}[theo]{Remark}
\newtheorem{prop}[theo]{Proposition}
\newtheorem{cor}[theo]{Corollary}
\newtheorem{ex}[theo]{Example}
\newtheorem{conj}[theo]{Conjecture}

\newcommand{\C}{\mathbb{C}}
\newcommand{\Qv}{Q^\vee}

\newcommand{\Z}{\mathbb{Z}}

\newcommand{\R}{\mathbb{R}}
\newcommand{\af}{\mathop{\mathrm{af}}}
\newcommand{\re}{\mathop{\mathrm{re}}}
\newcommand{\cl}{\mathrm{cl}}

\newcommand{\Waf}{W_{\af}}

\newcommand{\repo}{\Phi^{\re}_{\af, +}}

\newcommand{\qu}{\quad}
\newcommand{\si}{\frac{\infty}{2}}
\newcommand{\rp}{\mathcal{R}}
\newcommand{\inv}{^{-1}}
\newcommand{\F}{\mathcal{F}}
\newcommand{\la}{\langle}
\newcommand{\ra}{\rangle}
\newcommand{\h}{\mathcal{H}}
\newcommand{\haf}{\h_{\af}}
\newcommand{\nb}{\Delta}
\newcommand{\DBG}{\mathrm{DBG}}
\newcommand{\DBP}{\mathrm{DBP}}
\newcommand{\g}{\mathfrak{g}}
\newcommand{\ha}{\mathfrak{h}}
\newcommand{\wt}{\mathrm{wt}}
\newcommand{\sib}{<_{\si}}
\newcommand{\tw}{\mathop{\mathrm{tw}}}

\title[A combinatorial formula expressing periodic $R$-polynomials]{\Large{A combinatorial formula expressing periodic $R$-polynomials}}
\author[S. Naito and H. Watanabe]{Satoshi Naito and Hideya Watanabe}

\address{(S. Naito) Department of Mathematics, Tokyo Institute of Technology, 2-12-1 Oh-okayama, Meguro-ku, Tokyo 152-8551, Japan}
\email{naito@math.titech.ac.jp}

\address{(H. Watanabe) Department of Mathematics, Tokyo Institute of Technology, 2-12-1 Oh-okayama, Meguro-ku, Tokyo 152-8551, Japan}
\email{watanabe.h.at@m.titech.ac.jp}

\date{}

\begin{document}
\begin{abstract}
In 1980, Lusztig introduced the periodic Kazhdan-Lusztig polynomials, which are conjectured to have important information about the characters of irreducible modules of a reductive group over a field of positive characteristic, and also about those of an affine Kac-Moody algebra at the critical level. The periodic Kazhdan-Lusztig polynomials can be computed by using another family of polynomials, called the periodic $R$-polynomials. In this paper, we prove a (closed) combinatorial formula expressing periodic $R$-polynomials in terms of the ``doubled" Bruhat graph associated to a finite Weyl group and a finite root system.
\end{abstract}

\maketitle

\section{Introduction}
In \cite{KL}, Kazhdan and Lusztig introduced a family $\{P_{y,w}(q)\}_{y,w \in W}$ of polynomials, which are called the Kazhdan-Lusztig polynomials, for an arbitrary Coxeter system $(W,S)$. It is well-known that these polynomials have important representation-theoretic and geometric meanings. For example, it is known that the Kazhdan-Lusztig polynomials associated to a finite Weyl group describe the intersection cohomology of Schubert varieties. More precisely, let $G$ be a complex semisimple algebraic group, $B$ a Borel subgroup of $G$, and $W$ the associated (finite) Weyl group. Then the flag variety $X = G/B$ admits a decomposition (Bruhat decomposition)
\begin{align}
X = \bigsqcup_{w \in W} BwB/B. \nonumber
\end{align}
For each $w$, let $X_w$ denote the closure (with respect to the Zariski topology) of $C_w := BwB/B$ in $X$; the $X_w$'s are called the Schubert varieties. Then the Kazhdan-Lusztig polynomial $P_{y,w}(q)$ for $y,w \in W$ is identical to the Poincar\'{e} polynomial
\begin{align}
\sum_{i \geq 0} \dim \mathrm{IH}^{2i}_{C_y}(X_w)\  q^i; \nonumber
\end{align}
here, $\mathrm{IH}^{2i}(X_w)$ denotes the $2i$-th intersection cohomology sheaf on $X_w$, and $\mathrm{IH}^{2i}_{C_y}(X_w)$ is its stalk at a point of $C_y$. This equation also shows that the Kazhdan-Lusztig polynomials associated to an arbitrary finite Weyl group have no negative coefficients.

About the representation-theoretic side, we know that the Kazhdan-Lusztig polynomials describe the multiplicities of Verma modules in simple modules. Namely, let $\g$ be a finite-dimensional complex semisimple Lie algebra, $W$ its (finite) Weyl group, $P$ the weight lattice of $\g$, and $P_+$ the set of dominant integral weights. It is known that each finite-dimensional simple $\g$-module is isomorphic to the unique simple quotient $L(\mu)$ of the Verma module $V(\mu)$ for some $\mu \in P$. Kazhdan and Lusztig conjectured that for $\lambda \in P_+$ and $w \in W$, there holds the equality
\begin{align}
\mathrm{ch}L(w \circ \lambda) = \sum_{y \in W}(-1)^{\ell(y,w)}P_{y,w}(1) \ \mathrm{ch}M(y \circ \lambda), \nonumber
\end{align}
where $\ell:W \rightarrow \Z_{\geq0}$ denotes the length function on $W$, $\ell(y,w) := \ell(w) - \ell(y)$, and $W$ acts on $P$ by the so-called dot (or shifted) action. This conjecture was proved independently by Beilinson-Bernstein \cite{bebe} and Brylinski-Kashiwara \cite{bk}.

In the proof of the existence of the Kazhdan-Lusztig polynomials, another family of polynomials $\{R_{y,w}(q)\}_{y,w \in W}$, which we call the $R$-polynomials, play a key role. The Kazhdan-Lusztig polynomials and the $R$-polynomials are related by the following equation:
\begin{align}
q^{\ell(y,w)}\overline{P}_{y,w}(q) - P_{y,w}(q) = R_{y,w}(q) + \sum_{y<z<w} R_{y,z}(q)P_{z,w}(q)\ \text{for}\ y,w \in W. \label{eq:1}
\end{align}
Here, $\leq$ denotes the Bruhat order on $W$, and $\overline{f}(q) = f(q\inv)$ for $f(q) \in \Z[q,q\inv]$. Therefore, determining the Kazhdan-Lusztig polynomials is equivalent to determining the $R$-polynomials. As for the $R$-polynomials, Dyer proved the following combinatorial formula (in the notation of \cite[Chap. ${5}$]{bb}).
\begin{theo}[{\cite[(3.4)\ Corollary]{d},\ \cite[Theorem\ 5.3.4.]{bb}}]\label{bb}
Let $\preceq$ be a refection order on the set $\Phi_+$ of positive roots of an arbitrary finite Weyl group $W$. For $y,w \in W$, let $B^\prec (y,w)$ be the set of all sequences $\Delta = (y,ys_{\beta_1},ys_{\beta_1}s_{\beta_2},\ldots,ys_{\beta_1}s_{\beta_2}\cdots s_{\beta_{\ell(\Delta)}} = w)$ of elements of $W$ such that $y<ys_{\beta_1}<\ldots<ys_{\beta_1}\cdots ys_{\beta_{\ell(\Delta)}} = w$ and $\beta_1 \prec \cdots \prec \beta_{\ell(\Delta)}$, where for $\beta \in \Phi_+$, $s_\beta$ denotes the reflection with respect to $\beta \in \Phi_+$. Then, we have
\begin{align}
R_{y,w}(q) = \sum_{\nb \in B^\prec (y,w)} q^{\frac{1}{2}\bigl(\ell\left(y,w\right)-\ell\left(\nb\right)\bigr)}(q-1)^{\ell(\nb)}. \nonumber
\end{align}
\end{theo}

\begin{rem}\normalfont
Dyer proved Theorem \ref{bb} in a more general setting of an arbitrary Coxeter group.
\end{rem}

Also, in \cite{l1}, Lusztig introduced a family $\{ Q_{y,w}(q)\}_{y,w \in \Waf}$ of polynomials; here, $\Waf$ is the affine Weyl group. We call the $Q_{y,w}(q)$'s the periodic Kazhdan-Lusztig polynomials. These polynomials are also conjectured to have some representation-theoretic meaning. Let $G$ be a reductive group over a field of characteristic $p \gg h$, where $h$ denotes the Coxeter number for the associated root system. By fixing a maximal torus of $G$, we denote its character group by $X$, and set
\begin{align}
X_+ &:= \{ \lambda \in X \mid \lambda(\alpha^\vee) \geq 0\ \mathrm{for\ all\ simple\ coroots}\ \alpha^\vee \}\ (\mathrm{the\ set\ of\ dominant\ weights}), \nonumber\\
X_0 &:= \{ \lambda \in X \mid 0 < \lambda(\alpha^\vee) < p\ \mathrm{for\ all\ simple\ coroots}\ \alpha^\vee\}. \nonumber
\end{align}
Then, Lusztig's conjecture is the following:
\begin{conj}
Let $\lambda \in X_+$ and $w \in \Waf$ be such that $w \circ \lambda \in X_0$. Then, we have
\begin{align}
\mathrm{ch}L(w \circ \lambda) = \sum_{x} (-1)^{\ell(x,w)}P_{x,w}(1) \ \mathrm{ch}V(x \circ \lambda), \nonumber
\end{align}
where $V(x \circ \lambda)$ denotes the Weyl module of $G$ with simple head $L(x \circ \lambda)$, and the sum is over all $x \in \Waf$ such that $x \leq w$ and $x \circ \lambda \in X_+$.
\end{conj}
This conjecture was settled in \cite{AJS} for a sufficiently large $p \gg h$. In addition, it is proved in \cite{K} that Lusztig's conjecture is equivalent to the equality:
\begin{align}
\mathrm{ch}L(w \circ \lambda) = \sum_{x} (-1)^{\ell(x,w)}Q_{x,w}(1) \ \mathrm{ch}V(x \circ \lambda), \nonumber
\end{align}
where the sum is over the same $x$ as above.

There is another conjecture concerning periodic Kazhdan-Lusztig polynomials, called the Feigin-Frenkel conjecture. Let $\g_{\af}$ be an affine Kac-Moody algebra, $\ha_{\af}$ its Cartan subalgebra, and $\Waf$ the associated affine Weyl group. For each $\lambda \in \ha^*_{\af}$, let $\Delta(\lambda)$ denote the Verma module with highest weight $\lambda$ over $\g_{\af}$. Then, the restricted Verma module $\overline{\Delta}(\lambda)$ is defined to be the quotient module of $\Delta(\lambda)$ by a certain ideal of the Feigin-Frenkel center (see \cite[\S 1.4]{A} for details). Also, let $L(\lambda)$ denote the irreducible quotient of $\Delta(\lambda)$. Now, the Feigin-Frenkel conjecture is the following:
\begin{conj}
Let $\lambda \in \ha^*_{\af}$ be a weight at the critical level, i.e., $\lambda(K) = -\rho(K)$, where $K$ is the central element of $\g_{\af}$ and $\rho \in \ha_{\af}^*$ is the Weyl vector, and let $y,w \in \Waf$. We also assume that $\lambda \in \ha^*_{\af}$ is regular and dominant in the sense of \cite[\S 1.5]{A}. Then,
\begin{align}
[\overline{\Delta}(w \circ \lambda) : L(y \circ \lambda)] = Q_{y,w}(1), \nonumber
\end{align}
where $[\overline{\Delta}(w \circ \lambda) : L(y \circ \lambda)]$ denotes the Jordan-H\"{o}lder multiplicity of $L(y \circ \lambda)$ in $\overline{\Delta}(w \circ \lambda)$.
\end{conj}

When one considers the periodic Kazhdan-Lusztig polynomials, the role of the $R$-polynomials is replaced by the periodic $R$-polynomials $\{ \rp_{y,w}(q)\}_{y,w \in \Waf}$. Namely, Lusztig (\cite[(7.3.1)]{l1}) proved that for $y,w \in W_{\af}$,
\begin{align}
q^{\ell^\si(y,w)}\overline{Q}_{y,w}(q) - Q_{y,w}(q) = (-q)^{\ell^\si(y,w)}\overline{\rp}_{y,w}(q) + \sum_{y\sib z \sib w}(-q)^{\ell^\si(y,z)}\overline{\rp}_{y,z}(q)Q_{z,w}(q); \label{eq:2}
\end{align}
here, $\ell^\si$ denotes the semi-infinite length, $\ell^\si(y,w) := \ell^\si(w) - \ell^\si(y)$, and $\leq_\si$ the semi-infinite Bruhat order (defined later). Roughly speaking, arguments in the periodic case are parallel to those in the ordinary case if one replaces the ordinary length and Bruhat order by the semi-infinite length and semi-infinite Bruhat order, respectively. From this point of view, equation \eqref{eq:2} seems to be different from the equation obtained from equation \eqref{eq:1} by replacing $P$, $R$, $\ell$, and $\leq$ by $Q$, $\mathcal{R}$, $\ell^\si$, and $\leq_\si$, respectively. However, these two actually agree because
\begin{align}
R_{y,w}(q) = (-q)^{\ell(y,w)}\overline{R}_{y,w}(q). \nonumber
\end{align}
for $y,w \in W_{\af}$ (see, for example, \cite[Chap. 5, Exerc. 1]{bb}).

The periodic $R$-polynomials can also be thought of as a generalization of the ordinary $R$-polynomials in the following sense: if $y,w \in W$ (and if we regard the finite Weyl group $W$ as a subgroup of $\Waf$), then
\begin{align}
\rp_{y,w}(q) = R_{y,w}(q). \nonumber
\end{align}
Hence the periodic Kazhdan-Lusztig polynomials can be regarded as a generalization of the ordinary Kazhdan-Lusztig polynomials in the following sense: if $y,w \in W$, then
\begin{align}
Q_{y,w}(q) = P_{y,w}(q). \nonumber
\end{align}

The purpose of this paper is to give a (closed) combinatorial formula expressing periodic $R$-polynomials, which generalizes Theorem \ref{bb} in the following sense: if one applies our formula to elements of a finite Weyl group, then it agrees with the one in Theorem \ref{bb}. Namely, we prove the following theorem:
\begin{theo}\label{1.0.5}
Let $\preceq$ be a reflection order on $\Phi_+$, and $y,w \in \Waf$. Then we have
\begin{align}
\rp_{y,w}(q) = \sum_{\nb \in P^\prec (y,w)} q^{\deg(\nb)}(q-1)^{\ell(\nb)}; \nonumber
\end{align}
here, $P^\prec (y,w)$ consists of sequences $\Delta = (y = y_0 \xrightarrow[\beta_1]{} y_1 \xrightarrow[\beta_2]{} \cdots \xrightarrow[\beta_{\ell(\Delta)}]{} y_{\ell(\Delta)} = w)$ such that $y = y_0 \sib \cdots \sib y_{\ell(\Delta)} = w$, $\beta_1 \prec \cdots \prec \beta_{\ell(\Delta)}$, and some additional conditions $($see Definition \ref{2.4.2} for details$)$.
\end{theo}
If $y,w \in W$, then $P^\prec (y,w)$ agrees with $B^\prec (y,w)$ (in Theorem \ref{bb}), and for each $\nb \in P^\prec (y,w) = B^\prec (y,w)$, we obtain $\deg(\nb) = \frac{1}{2}(\ell(y,w) - \ell(\nb))$. We prove Theorem \ref{1.0.5} by constructing explicit bijections between (specific subsets of) $P^\prec (y,w)$, $y,w \in W_{\af}$. In particular, we construct explicit bijections $B^\prec (sy,sw) \sqcup B^\prec (sy,w) \rightarrow B^\prec (y,w)$ for all $y,w \in W,\ s \in S$ such that $y < sy,\ sw < w$, and bijections $B^\prec (sy,sw) \rightarrow B^\prec (y,w)$ for all $y,w \in W,\ s \in S$ such that $sy < y,\ sw < w$, where $S$ denotes the set of simple reflections for $W$. We remark that this also gives a new proof of Theorem \ref{bb}.

Now, let us recall that Dyer introduced the notion of twisted Bruhat order of a Coxeter symstem in \cite{d1}. He also defined polynomials $R_{\leq_{\tw}}(x,y)$ for a twisted Bruhat order $\leq_{\tw}$ and $x, y \in W$ such that each length two closed subinterval of $[x,y]$ with respect to $\leq_{\tw}$ has exactly four elements. In \cite{d3}, it is shown that the semi-infinite Bruhat order is a twisted Bruht order. Then, we remark that the periodic $R$-polynomials $\mathcal{R}_{x,y}$ coincide with Dyer's polynomials $R_{\leq_\si}(x, y)$ for all $x,y \in \Waf$ such that each length two closed subinterval of $[x,y]$ with respect to $\leq_\si$ has exactly four elements. This is because if $x, y \in \Waf$ satisfy the condition on the  lenght two closed subintervals, then both polynomials have the same recursive formulas $(R1) - (R3)$ in Proposition \ref{periodicR}, and are uniquely determined by them. With our terminology (in Definition \ref{2.4.2}), the condition on the subintervals is equivalent to that each $\Delta \in P^\prec (z,w)$ has no translation edges for all $z, w \in [x,y]$ (see Proposition \ref{atom}). This condition fails if there are $z, w \in [x,y]$ such that $\cl(z) = \cl(w)$. In particular, if $\ell^\si(x,y) > |W|$, then there exist $z,w \in [x,y]$ such that $\cl(z) = \cl(w)$, and therefore periodic $\mathcal{R}$-polynomial $\mathcal{R}_{x,y}$ cannot be written as Dyer's $R_{\leq_{\tw}}(x,y)$.

This paper is organized as follows. In Section $2$, we fix our notation and recall the notions and basic properties of affine root systems, reflection orders, semi-infinite Bruhat order, $R$-polynomials, and periodic $R$-polynomials. Also, we state our main result (Theorem \ref{main}). Section $3$ is devoted to the preparation for the proof of Theorem \ref{main}. In Section $4$, we complete the proof of Theorem \ref{main}. In Section $5$, we introduce a directed graph, which we call the double Bruhat graph, in order to describe Theorem \ref{main} in a way more suitable for actual computation.

\section{Setting and the main result}
In this section, we fix our notation and state our main result (Theorem \ref{main}).
\subsection{Affine root systems and affine Weyl groups}
Let $\g$ be a finite-dimensional, complex simple Lie algebra of rank $l$, and $\ha$ a Cartan subalgebra of $\g$. We denote by $\Phi \subset \ha^*$ the set of roots of $\g$ with respect to $\ha$, by $\Pi = \{ \alpha_1,\ldots,\alpha_l \}$ a set of simple roots, by $\Phi_+$ the set of positive roots with respect to $\Pi$, and by $\{ \alpha_1^\vee,\ldots,\alpha_l^\vee \} \subset \ha$ the set of simple coroots. Also, we denote by $\la\ ,\ \ra:\ha^* \times \ha \rightarrow \C$ the duality pairing. For each $\alpha \in \Phi$, let $\alpha^\vee \in \ha$ denote the coroot of $\alpha$. Let $\Qv = \sum_{i=1}^l \Z\alpha_i^\vee$ be the coroot lattice, with $\Qv_+ = \sum_{i=1}^l \Z_{\geq 0}\alpha_i^\vee$, and let $P^\vee = \{ \mu \in \ha \mid \la \alpha_i,\mu \ra \in \Z \ \mathrm{for\ all}\ i = 1,\ldots,l \}$ be the coweight lattice, with $P^\vee_+ = \{ \mu \in P^\vee \mid \la \alpha_i,\mu \ra \geq 0 \ \mathrm{for\ all}\ i=1,\ldots,l \}$ the set of dominant coweights. 

Let $W$ denote the (finite) Weyl group of $\g$, with $S = \{ s_1,\ldots,s_l \}$ the set of simple reflections. Then $W$ and $\Qv$ act on $\ha$ on the left: for each $i = 1,\ldots,l$, $\lambda \in \Qv$, and $\mu \in \ha$,
\begin{align}
s_i \cdot \mu = \mu - \la \alpha_i,\mu \ra \alpha_i^\vee, \qu t_\lambda \cdot \mu = \mu + \lambda, \nonumber
\end{align}
where $t_\lambda$ is the image of $\lambda \in \Qv$ in the group $\mathrm{Aut}(\ha)$ of all affine transformations on $\ha$. Since these actions are faithful, we can regard $W$ and $\Qv$ as subgroups of $\mathrm{Aut}(\ha)$. Also, let $\theta \in \Phi_+$ denote the highest root, and set $s_0 = s_{\theta}t_{-\theta^\vee} \in \mathrm{Aut}(\ha)$, where $s_\theta$ is the reflection with respect to $\theta$. Then, $S_{\af} := \{s_0,s_1,\ldots,s_l \}$ generate the affine Weyl group $\Waf \simeq W \ltimes \Qv$.

Now, let $\g_{\af} := \g \otimes \C[t,t\inv] \oplus \C K \oplus \C D$ be the affine Lie algebra associated to $\g$, where $K$ is the canonical central element; the corresponding Cartan subagebra is
\begin{align}
\ha_{\af} := \ha \oplus \C K \oplus \C D. \nonumber
\end{align}
We extend the duality pairing $\la \ ,\ \ra : \ha^* \times \ha \rightarrow \C$ to the duality pairing $\la \ , \ \ra : \ha_{\af}^* \times \ha_{\af} \rightarrow \C$ in such a way that $\la \ha^* , \C K \oplus \C D \ra = 0$. Let $\delta \in \ha_{\af}^*$ be the element defined by
\begin{align}
\la \delta , \ha \oplus \C K \ra = 0, \nonumber\\
\la \delta , D \ra = 1. \nonumber
\end{align}
Then, the set $\repo$ of positive real roots of $\g_{\af}$ can be written as $\repo = \{ \alpha + m\delta \mid \alpha \in \Phi_+,\ m \in \Z_{\geq 0} \} \sqcup \{ -\alpha + m\delta \mid \alpha \in \Phi_+,\ m \in \Z_{> 0} \}$, and the set $\Phi^{\re}_{\af}$ of real roots is $\repo \sqcup (-\repo)$. Also, the set of simple roots of $\g_{\af}$ is $\{ \alpha_0,\alpha_1,\ldots,\alpha_l \}$, where $\alpha_0 = -\theta + \delta$. Then, $\Waf$ acts on $\Phi^{\re}_{\af}$ by
\begin{align}
xt_\lambda (\alpha + n \delta) = x \alpha + (n - \la \lambda, \alpha \ra) \delta \nonumber
\end{align}
for $x \in W$, $\lambda \in Q^\vee$, $\alpha \in \Phi$, and $n \in \Z$. For each $\beta \in \Phi^{\re}_{\af}$, let $\beta^\vee \in \ha_{\af}$ denote the coroot of $\beta$.

For each $\alpha \in \Phi_+$ and $m \in \Z$, let $\F_{\alpha,m} \subset \ha_{\R} := \R \otimes_{\Z} Q^\vee$ be the hyperplane defined by
\begin{align}
\F_{\alpha,m} = \{ \mu \in \ha_{\R} \mid \la \alpha,\mu \ra = m\}. \nonumber
\end{align}
In this notation, the reflection $s_{\pm \alpha + m \delta} \in W_{\af}$ corresponding to a positive real root $\pm\alpha + m\delta \in \repo$ is the reflection fixing the hyperplane $\F_{\alpha,\pm m}$ pointwise. An alcove is a connected component of $\ha_{\R} \setminus (\bigcup_{\alpha \in \Phi_+, m \in \Z}\ \F_{\alpha,m})$. The set $A^- := \{ \lambda \in \ha_{\R} \mid -1 < \la \alpha,\lambda \ra < 0\ \mathrm{for\ all}\ \alpha \in \Phi_+ \}$ is an alcove. It is known that the action of $\Waf$ on $\ha$ defined above induces a transitive and faithful right action of $\Waf$ on the set of alcoves, that is, the map $w := \cl(w)t_{\wt(w)} \mapsto A^- \cdot w := \cl(w)\inv A^- + \wt(w)$ is a bijection between $\Waf$ and the set of alcoves, where $\cl(w) \in W,\ \wt(w) \in \Qv$. Through this bijection, we sometimes identify these two sets.

For $i = 1,\ldots,l$ (resp., $i = 0$), the $i$-wall of $A^-$ is defined to be the intersection of the hyperplane $\F_{\alpha_i,0}$ (resp., $\F_{\theta,-1}$) and the closure $\overline{A^-}$ of $A^-$. The $i$-wall of an alcove $A$ is the unique element of the orbit of the $i$-wall of $A^-$ under the action of $\Waf$ that is contained in the closure $\overline{A}$ of $A$.

Here, we recall some basic facts about (affine) Weyl groups that are frequently used in this paper.

\begin{prop}
\begin{enumerate}
\item[$(1)$] $\ell(w) = \sharp \{ \alpha \in \Phi_+ \mid w(\alpha) \in - \Phi_+ \}$.
\item[$(2)$] $($\cite{h2}$)$ Let $\alpha \in \Phi_+$ and $w \in W$. Then, $\ell(w) < \ell(s_\alpha w)$ if and only if $w\inv(\alpha) \in \Phi_+$.
\item[$(3)$] $($\cite{bfp}$)$ For each $\alpha \in \Phi_+$, we have $\ell(s_\alpha) \leq 2 \la \rho, \alpha^\vee \ra - 1$.
\end{enumerate}
\end{prop}

\subsection{The semi-infinite Bruhat order}
The finite Weyl group $W$ and the affine Weyl group $\Waf$ are equipped with a partial order $\leq$, which is known as the Bruhat order. Following \cite{l1}, we introduce another partial order $\leq_\si$ on $\Waf$, which we call the semi-infinite Bruhat order. First, for each $\alpha \in \Phi_+$ and $m \in \Z$, we divide $\ha_{\R}$ into three parts:
\begin{align}
\ha_{\R} &= \F^-_{\alpha,m} \sqcup \F_{\alpha,m} \sqcup \F^+_{\alpha,m}, \nonumber
\end{align}
where $\F^{\pm}_{\alpha,m} = \{ \mu \in \ha_{\R} \mid \pm(\la \alpha,\mu \ra - m) > 0 \}$. Next, to each pair $(A,B)$ of alcoves, we associate an integer $d(A,B) \in \Z$ as follows. Let $(A=A_0,A_1,\ldots,A_r=B)$ be a sequence of alcoves such that $A_{i-1}$ and $A_i$ have a unique common wall for each $i$; there exists a unique $\beta_i \in \Phi_+$ and $m_i \in \Z$ such that $\F_{\beta_i,m_i}$ contains the common wall. Then we define $d(A,B)$ by:
\begin{align}
d(A,B) &= \sum_{i=1}^{r} d(A_{i-1},A_i), \nonumber\\
d(A_{i-1},A_i) &= \begin{cases}
1 \qu & \mathrm{if}\ A_i \subset \F^+_{\beta_i,m_i}\ (\mathrm{or\ equivalently},\ A_{i-1} \subset \F^-_{\beta_i,m_i}),\\
-1 \qu & \mathrm{if}\ A_i \subset \F^-_{\beta_i,m_i}\ (\mathrm{or\ equivalently},\ A_{i-1} \subset \F^+_{\beta_i,m_i}).
\end{cases} \nonumber
\end{align}
For $y,w \in \Waf$, we set $\ell^\si(y,w) := d(A^- \cdot y,A^- \cdot w)$. In particular, when $y$ is the identity element $e$, we simply write $\ell^\si(e,w) = \ell^\si(w)$ and call it the semi-infinite length of $w$. By definition, we have $\ell^\si(y,w) = \ell^\si(w) - \ell^\si(y)$.

\begin{ex}\normalfont
In type $A_2$:\\
\begin{figure}[H]
\begin{tabular}{cc}
\begin{minipage}{0.5\hsize}
\begin{center}
\labellist
\small\hair 2pt
\pinlabel $e$ at 180 180
\pinlabel $s_1$ at 260 180
\pinlabel $s_2$ at 140 245
\pinlabel $s_2s_1$ at 285 245
\pinlabel $s_1s_2$ at 175 305
\pinlabel $w_0$ at 260 305
\pinlabel $s_0$ at 145 115
\pinlabel $t_{\alpha_1^\vee}$ at 400 180
\pinlabel $t_{\alpha_2^\vee}$ at 65 370
\pinlabel $t_{\theta^\vee}$ at 290 370
\endlabellist
$\includegraphics[width=5cm]{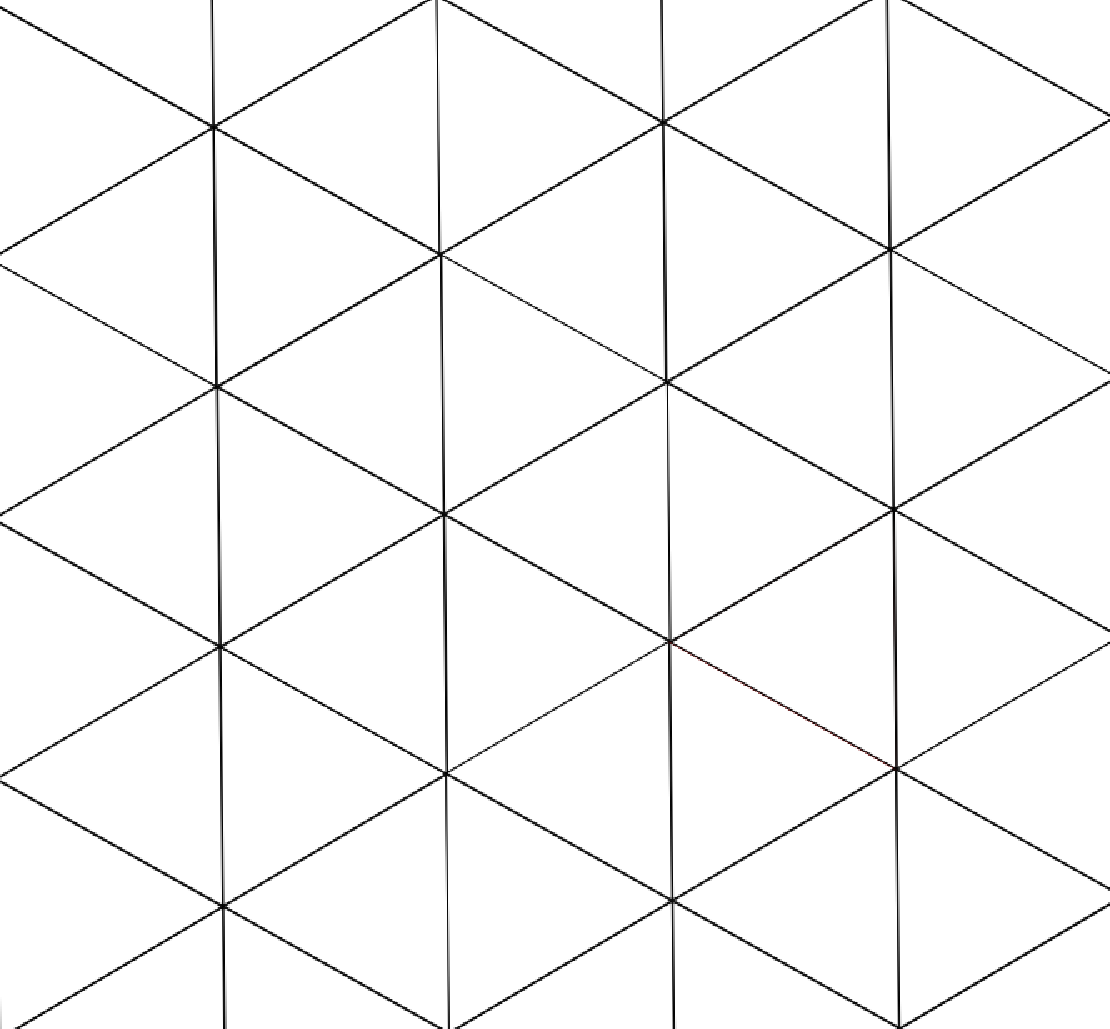}$
\caption{alcoves and elements of $\Waf$}
\label{fig:alcoves}
\end{center}
\end{minipage}
\begin{minipage}{0.5\hsize}
\begin{center}
\labellist
\small\hair 2pt
\pinlabel $0$ at 180 180
\pinlabel $1$ at 260 180
\pinlabel $1$ at 140 245
\pinlabel $2$ at 285 245
\pinlabel $2$ at 175 305
\pinlabel $3$ at 260 305
\pinlabel $-1$ at 145 115
\pinlabel $2$ at 390 180
\pinlabel $2$ at 65 370
\pinlabel $4$ at 275 370
\endlabellist
$\includegraphics[width=5cm]{alcoves2}$
\caption{semi-infinite lengths}
\label{fig:silength}
\end{center}
\end{minipage}
\end{tabular}
\end{figure}
\end{ex}

Finally, we define the partial order $\leq_\si$ on $\Waf$ as follows: for $y,w \in \Waf$, we write $y \leq_\si w$ if there exists a sequence $(y=y_0,y_1,\ldots,y_r=w)$ of elements of $\Waf$ such that $y_i = y_{i-1}s_{\beta_i}$ for some $\beta_i \in \repo$ and $\ell^\si(y_{i-1},y_i) > 0$ for each $1 \leq i \leq r$. Note that if $y \sib w$, then $\ell^\si(y) < \ell^\si(w)$.

Let $\mathcal{C}^+$ denote the dominant Weyl chamber:
\begin{align}
\mathcal{C}^+ = \{ \mu \in \ha_{\R} \mid \la \alpha_i,\mu \ra > 0 \ \mathrm{for\ all}\ i = 1,\ldots,l \}. \nonumber
\end{align}

Let $\rho$ denote half the sum of positive roots of $\g$. Then,  $\rho$ is characterized by the condition that $\la \rho, \alpha_i^\vee \ra = 1$ for all $i \in I$.

About the semi-infinite Bruhat order and the semi-infinite length, we know the following.

\begin{prop}[{\cite[Lecture 13]{P}}]\label{peterson}
\begin{enumerate}
\item[$(1)$] Let $w = \cl(w)t_{\wt(w)} \in \Waf$ $(\cl(w) \in W,\ \wt(w) \in \Qv)$. Then $\ell^\si(w) = \ell(\cl(w)) + 2\la \rho, \wt(w) \ra$.
\item[$(2)$] Let $y,w \in \Waf$. Then $y \leq_\si w$ if and only if $yt_\lambda \leq wt_\lambda$ for some $\lambda \in \Qv_+$ such that $A^- \cdot yt_\lambda, A^- \cdot wt_\lambda \subset \mathcal{C}^+$.
\end{enumerate}
\end{prop}

\begin{prop}\label{atom}
Let $y, w \in \Waf$ be such that $\ell^\si(y,w) = 2$. Then the number of elements $z \in \Waf$ satisfying $y <_\si z <_\si w$ is $1$ if $\cl(y) = \cl(w)$, and is $2$ if $\cl(y) \neq \cl(w)$.
\end{prop}

\begin{proof}
By Proposition \ref{peterson} (1), we have $\ell^\si(y,ys_{\alpha + n\delta})$ is odd for all $\alpha + n\delta \in \Phi^{\re}_{\af, +}$. Since $\ell^\si(y,w) = 2$, by the definition of $\leq_\si$, there exist $z \in \Waf$ and $\alpha + n\delta, \beta + m\delta \in \Phi^{\re}_{\af,+}$ such that $z = y s_{\alpha + n\delta}$, $w = z s_{\beta + m\delta}$, and $\ell^\si(y,z) = \ell^\si(z,w) = 1$.

Assume that $\cl(y) = \cl(w)$, equivalently $\alpha = \beta$. In this case, we have $w = y t_{l\alpha}$ for some $l \in \Z$. Then, $\ell^\si(y,w) = 2 \la \rho, l \alpha \ra$. Since $\ell^\si(y,w) = 2$, we have $l = 1$ and $\alpha$ is a simple root, say $\alpha_i$ $(i \in I)$. This implies $w = y t_{\alpha_i + n\delta}$ and $z = y s_{\alpha_i + n\delta}$. In particular, $z$ is uniquely determined by $y$ and $w$. Therefore, the number of elements $z \in \Waf$ satisfying $y <_\si z <_\si w$ is exactly $1$.

Next, assume that $\cl(y) \neq \cl(w)$. By {\cite[Proposition A.1.2]{ins}} , the number in question is equal to the number of paths of length two from $\cl(y)$ to $\cl(w)$ in the directed graph $D(W)$ of {\cite[Definition 6.2]{bfp}}. Each path of length two from $\cl(y)$ to $\cl(w)$ in $D(W)$ is either label increasing or label decreasing. Then, by {\cite[Theorem 6.6.1]{bfp}}, the number of such paths is exactly $2$, which completes the proof of the proposition.
\end{proof}

\begin{lem}[{\cite[Lecture 13, Proposition 1]{P}}]\label{223}
Let $w \in \Waf$ and $\mu \in \Qv_+$. Then we have
\begin{align}
w \leq_\si wt_\mu. \nonumber
\end{align}
\end{lem}

\begin{proof}
It suffices to show the assertion when $\mu$ is a simple coroot, that is, when $\mu = \alpha_i^\vee$ for some $i = 1,\ldots,l$. We write $w = \cl(w)t_{\wt(w)}$. Assume first that $\cl(w) < \cl(w)s_i$. Then, we have
\begin{align}
\cl(w)s_it_{\wt(w)} = \cl(w)t_{s_i(\wt(w))}s_i &= \cl(w)t_{\wt(w)- \la \alpha_i,\wt(w)\ra \alpha_i^\vee}s_i \nonumber\\
&= \cl(w)t_{\wt(w)}s_it_{\la \alpha_i,\wt(w)\ra \alpha_i^\vee} = ws_{\alpha_i + \la \alpha_i,\wt(w) \ra \delta}, \nonumber
\end{align}
and
\begin{align}
\ell^\si(\cl(w)s_it_{\wt(w)}) = \ell(\cl(w)s_i) + 2\la \rho,\wt(w) \ra = \ell(\cl(w)) + 1 + 2\la \rho,\wt(w) \ra = \ell^\si(w) + 1. \nonumber
\end{align}
Therefore, we obtain $w \sib \cl(w)s_it_{\wt(w)}$. Furthermore, we have
\begin{align}
wt_{\alpha_i^\vee} = \cl(w)s_is_it_{\wt(w)+\alpha_i^\vee} = \cl(w)s_it_{\wt(w)}s_it_{(\la \alpha_i,\wt(w) \ra +1)\alpha_i^\vee} = \cl(w)s_it_{\wt(w)}s_{\alpha_i + (\la \alpha_i,\wt(w) \ra +1)\delta}, \nonumber
\end{align}
and
\begin{align}
\ell^\si(wt_{\alpha_i^\vee}) = \ell(\cl(w)) + 2\la \rho,\wt(w)+\alpha_i^\vee \ra &= \ell(\cl(w)) + 2\la \rho,\wt(w)\ra + 2 \nonumber\\
&= \ell^\si(w) + 2 = \ell^\si(\cl(w)s_it_{\wt(w)}) + 1. \nonumber
\end{align}
From these, we deduce that $\cl(w)s_it_{\wt(w)} \sib wt_{\alpha_i^\vee}$, and hence
\begin{align}
w \sib wt_{\alpha_i^\vee}. \nonumber
\end{align}
Also, we obtain $\cl(w)s_it_{\wt(w)} \sib \cl(w)s_it_{\wt(w)+\alpha_i^\vee}$, since $wt_{\alpha_i^\vee} = \cl(w)t_{\wt(w)+\alpha_i^\vee} \sib \cl(w)s_it_{\wt(w) + \alpha_i^\vee}$. This implies the assertion when $\cl(w)s_i < \cl(w)$. This proves the lemma.
\end{proof}

\begin{prop}[{\cite[Lecture 13, Proposition 2]{P}}]\label{criterion}
Let $w = \cl(w)t_{\wt(w)}\in \Waf$, and $\beta = \alpha + m\delta \in \repo$. Then, $\ell^\si(w) < \ell^\si(s_\beta w)$ if and only if $\cl(w)\inv(\alpha) \in \Phi_+$.
\end{prop}

\begin{proof}
By direct calculation, we see that
\begin{align}
s_\beta w &= \cl(w)s_{\cl(w)\inv(\beta)}t_{\wt(w)} \nonumber\\
&= \cl(w)s_{\cl(w)\inv(\alpha)}t_{\wt(w) + m\cl(w)\inv(\alpha^\vee)}. \nonumber
\end{align}
Hence it follows that
\begin{align}
\ell^\si(s_\beta w) - \ell^\si(w) &= \ell(\cl(w)s_{\cl(w)\inv(\alpha)}) + 2\la \rho,\wt(w) + m\cl(w)\inv(\alpha^\vee) \ra - (\ell(\cl(w)) + 2\la \rho,\wt(w) \ra) \nonumber\\
&= \ell(\cl(w)s_{\cl(w)\inv(\alpha)}) - \ell(\cl(w)) + 2m \la \rho,\cl(w)\inv(\alpha^\vee) \ra. \nonumber
\end{align}

We first prove the ``only if" part. Assume that $\ell^\si(w) < \ell^\si(s_\beta w)$ and $\cl(w)\inv(\alpha) \in -\Phi_+$. Then $\ell(s_{\cl(w)\inv(\alpha)}) \leq -2 \la \rho,\cl(w)\inv(\alpha^\vee) \ra -1$, and
\begin{align}
0 < \ell^\si(s_\beta w) - \ell^\si(w) &= \ell(\cl(w)s_{\cl(w)\inv(\alpha)}) - \ell(\cl(w)) + 2m\la\rho,\cl(w)\inv(\alpha^\vee)\ra \nonumber\\
&\leq \ell(\cl(w)) + \ell(s_{\cl(w)\inv(\alpha)}) - \ell(\cl(w)) + 2m\la \rho,\cl(w)\inv(\alpha^\vee)\ra \nonumber\\
&\leq -2\la \rho,\cl(w)\inv(\alpha^\vee)\ra -1 + 2m\la\rho,\cl(w)\inv(\alpha^\vee)\ra \nonumber\\
&< 2(m-1)\la\rho,\cl(w)\inv(\alpha^\vee)\ra. \nonumber
\end{align}
This implies that $m = 0$ since $\la \rho,\cl(w)\inv(\alpha^\vee)\ra < 0$. Therefore, $\alpha = \beta \in \repo$, and hence $\alpha \in \Phi_+$. Hence it follows that
\begin{align}
\ell^\si(s_\beta w) - \ell^\si(w) &= \ell(s_\alpha \cl(w)) - \ell(\cl(w)) < 0 \qu (\mathrm{since}\  \alpha \in \Phi_+\ \mathrm{and}\ \cl(w)\inv(\alpha) \in -\Phi_+), \nonumber
\end{align}
which is a contradiction.

Next we prove the ``if" part. Assume that $\cl(w)\inv(\alpha) \in \Phi_+$. From the above, we have
\begin{align}
\ell^\si(s_\beta w) - \ell^\si(w) = \ell(\cl(w)s_{\cl(w)\inv(\alpha)}) - \ell(\cl(w)) + 2m\la\rho,\cl(w)\inv(\alpha^\vee)\ra; \nonumber
\end{align}
note that $\la\rho,\cl(w)\inv(\alpha^\vee)\ra > 0$. If $\alpha \in \Phi_+$, then $\cl(w)s_{\cl(w)\inv(\alpha)} = s_\alpha \cl(w) > \cl(w)$, and hence $\ell^\si(s_\beta w) - \ell^\si(w) > 0$. If $\alpha \in -\Phi_+$, then $m > 0$ and hence
\begin{align}
\ell^\si(s_\beta w) - \ell^\si(w) &= \ell(\cl(w)s_{\cl(w)\inv(\alpha)}) - \ell(\cl(w)) + 2m\la\rho,\cl(w)\inv(\alpha^\vee)\ra \nonumber\\
&\geq \ell(\cl(w)) - \ell(s_{\cl(w)\inv(\alpha)}) - \ell(\cl(w)) + 2m\la\rho,\cl(w)\inv(\alpha^\vee)\ra \nonumber\\
&\geq -2\la\rho,\cl(w)\inv(\alpha^\vee)\ra +1 + 2m\la\rho,\cl(w)\inv(\alpha^\vee)\ra \nonumber\\
&= 2(m-1)\la\rho,\cl(w)\inv(\alpha^\vee)\ra + 1 > 0, \nonumber
\end{align}
as desired. This proves the lemma.
\end{proof}

\subsection{$R$-polynomials and periodic $R$-polynomials}
Let $\h$ be the Hecke algebra of $W$, and $\haf$ the Hecke algebra of $\Waf$, that is, $\h$ (resp., $\haf$) is the associative algebra over $\Z[q,q\inv]$ generated by $\{ T_s \mid s \in S \}$ (resp., $\{ T_s \mid s \in S_{\af} \}$) subject to the braid relations and the Hecke relation
\begin{align}
T_s^2 = q + (q-1)T_s\ \mathrm{for}\ s \in S\ (\mathrm{resp}.,\ s \in S_{\af}). \nonumber
\end{align}
For each $x \in W$ (resp., $w \in \Waf$), we set $T_x = T_{s_{i_1}}\cdots T_{s_{i_k}}$ (resp., $T_w = T_{s_{i_1}} \cdots T_{s_{i_k}}$), where $s_{i_1}\cdots s_{i_k}$ is a reduced expression for $x$ (resp., $w$). Let $\bar{\ }:\h \rightarrow \h$ (resp., $\bar{\ }:\haf \rightarrow \haf$) denote the involutive $\Z$-algebra automorphism that sends $q$ to $q\inv$ and $T_s$ to $T_s\inv = q\inv T_s + (q\inv -1)$.

\begin{defi}\normalfont
Let $u,v \in W$. Then, the $R$-polynomial $R_{u,v}(q) \in \Z[q,q\inv]$ is defined by
\begin{align}
\overline{T}_{v} = q^{-\ell(v)}\sum_{u \in W} (-1)^{\ell(u,v)} R_{u,v}(q) T_u. \nonumber
\end{align}
\end{defi}
It follows from the next proposition that $R_{u,v}(q)$ is, in fact, an element of $\Z[q]$ of degree $\ell(u,v)$.

\begin{prop}[{\cite[5.3]{bb}, \cite[7.5]{h2}}]\label{recursive}
The $R$-polynomials are uniquely determined by the following equations.
\begin{enumerate}
\item[$(1)$] $R_{v,v}(q) = 1$ for all $v \in W$.
\item[$(2)$] $R_{u,v}(q) = 0$ unless $u \leq v$ $(u,v \in W)$.
\item[$(3)$] $R_{u,v}(q) = \begin{cases}
R_{su,sv}(q) \qu & \mathrm{if}\ sv<v\ \mathrm{and}\ su<u,\\
qR_{su,sv}(q) + (q-1)R_{u,sv}(q) \qu & \mathrm{if}\ sv<v\ \mathrm{and}\ u<su,
\end{cases}$\\
for all $u,v \in W$ and all $s \in S$.
\end{enumerate}
Also, equation $(3)$ can be replaced by
\begin{enumerate}
\item[$(3')$] $R_{u,v}(q) = \begin{cases}
R_{ut,vt}(q) \qu & \mathrm{if}\ vt<v\ \mathrm{and}\ ut<u,\\
qR_{ut,vt}(q) + (q-1)R_{u,vt}(q) \qu & \mathrm{if}\ vt<v\ \mathrm{and}\ u<ut,
\end{cases}$\\
for all $u,v \in W$ and all $t \in S$.
\end{enumerate}
\end{prop}


Now, following \cite{K} (see also \cite{l1}), we define the periodic $R$-polynomials. For each $\lambda \in \Qv$, there exists $\mu \in \Qv \cap P^{\vee}_+$ such that $\lambda + \mu \in \Qv \cap P^{\vee}_+$. Then we set
\begin{align}
X^\lambda := T_{t_{\lambda + \mu}} T_{t_\mu}\inv \in \h_{\af}; \nonumber
\end{align}
it is known that this definition of $X^\lambda$ is independent of the choice of $\mu$. Also, for each $w = \cl(w)t_{\wt(w)}$ $(\cl(w) \in W, \ \wt(w) \in \Qv)$, we set
\begin{align}
\widetilde{T}_w := T_{\cl(w)}X^{\wt(w)}. \nonumber
\end{align}
Then, $\{ \widetilde{T}_{w} \mid w \in \Waf \}$ forms a basis of $\haf$. With respect to this basis, we enlarge $\haf$ as follows. We set
\begin{align}
\widehat{\haf} := \left\{ \sum_{w \in \Waf} a_w \widetilde{T}_w \mid \{ w \in \Waf \mid a_w \neq 0\}\ \text{is bounded from above} \right\}, \nonumber
\end{align}
where a subset $Y \subset \Waf$ is said to be bounded from above if there exists $w \in \Waf$ such that $y \sib w$ for all $y \in Y$; note that $\widehat{\haf}$ is no longer an algebra, but is still a left $\haf$-module.

\begin{theo}[{\cite[Theorem 2.12]{l1}}, {\cite[Proposition 2.8]{K}}]\label{234}
There exists a unique involutive $\Z$-linear automorphism $\Psi:\widehat{\haf} \rightarrow \widehat{\haf}$ such that
\begin{align}
\Psi(h\cdot m) = \overline{h} \cdot \Psi(m) \qu (h \in \haf,\ m \in \widehat{\haf}), \nonumber
\end{align}
and
\begin{align}
\Psi \Bigl(\sum_{x \in W} \widetilde{T}_{xt_\lambda} \Bigr) = q^{-\ell^\si(w_0t_\lambda)}\sum_{x \in W} \widetilde{T}_{xt_\lambda} \qu (\lambda \in \Qv). \nonumber
\end{align}
Here, $w_0 \in W$ denotes the longest element.
\end{theo}

\begin{defi}\normalfont
Let $y,w \in \Waf$. Then the periodic $R$-polynomial $\rp_{y,w}(q)$ is defined by
\begin{align}
\Psi(\widetilde{T}_w) = q^{-\ell^\si(w)}\sum_{y \in \Waf} (-1)^{\ell^\si(y,w)} \rp_{y,w}(q) \widetilde{T}_y. \nonumber
\end{align}
\end{defi}
It follows from the next proposition that $\rp_{y,w}(q)$ is, in fact, an element of $\Z[q]$ of degree $\ell^\si(y,w)$.

\begin{prop}[{\cite[Proposition 11.1]{l1}}]\label{periodicR}
The periodic $R$-polynomials are uniquely determined by the following equations.
\begin{enumerate}
\item[$(\mathrm{R}1)$] $\rp_{w,w}(q) = 1$ for all $w \in \Waf$.
\item[$(\mathrm{R}2)$] $\rp_{y,w}(q) = 0$ unless $y \leq_\si w$ $(y,w \in \Waf)$.
\item[$(\mathrm{R}3)$] $\rp_{y,w}(q) = \begin{cases}
\rp_{sy,sw}(q) \qu & \mathrm{if}\ sw \sib w\ \mathrm{and}\ sy \sib y, \\
q\rp_{sy,sw}(q) + (q-1)\rp_{y,sw}(q) \qu & \mathrm{if}\ sw \sib w\ \mathrm{and}\ y \sib sy,
\end{cases}$\\
for all $y,w \in \Waf$ and all $s \in S_{\af}$.
\item[$(\mathrm{R}4)$] \begin{align}
\sum_{x \in W} (-1)^{\ell^\si(y,xt_\lambda)}q^{\ell(xw_0)}\rp_{y,xt_\lambda}(q) = \delta_{\lambda,\lambda(y)}\ \ \mathrm{for\ all}\ \lambda \in \Qv,\ y \in \Waf \nonumber.\end{align}
\end{enumerate}
\end{prop}

\begin{rem}\normalfont
Equation (R4) is equivalent to the second equation in Theorem $\ref{234}$:
\begin{align}
\Psi \Bigl(\sum_{x \in W} \widetilde{T}_{xt_\lambda} \Bigr) = q^{-\ell^\si(w_0t_\lambda)}\sum_{x \in W} \widetilde{T}_{xt_\lambda} \qu (\lambda \in \Qv). \nonumber
\end{align}
\end{rem}

\subsection{Main result}\label{sec;reforder}
In this subsection, we recall the notion of reflection orders and state our main result of this paper.
\begin{defi}\normalfont
A total order $\preceq$ on $\Phi_+$ $($resp., $\repo)$ is called a reflection order if it satisfies the following:
\begin{align}
\text{if}\ a,b \in \R_{>0},\alpha,\beta, a\alpha + b\beta \in \Phi_+\ (\text{resp}.,\ \repo)\ \text{and}\ \alpha \prec \beta,\ \text{then}\ \alpha < a\alpha + b\beta \prec \beta. \nonumber
\end{align}
\end{defi}

In the following, we take and fix a reflection order $\preceq$ on $\Phi_+$.

\begin{defi}\normalfont\label{2.4.2}
Let $y = \cl(y)t_{\wt(y)},w = \cl(w)t_{\wt(w)} \in \Waf$, and $\beta \in \Phi_+$.
\begin{enumerate}
\item[$(1)$] We write $y \xrightarrow[\beta]{m} w$ if either $w = yt_{m\beta^\vee}$ for some $m \in \Z_{>0}$, or $w = \cl(y)s_\beta t_{\wt(y) + m\beta^\vee}$ for some $m \in \Z_{\geq 0}$ and $y \sib w$. In the former case, $y$ and $w$ differ only by a translation, and in the latter case, $w$ is obtained from $y$ by multiplying the reflection $s_{\beta + (m + \la \beta,\wt(y) \ra)\delta}$ on the right, where $\beta + (m + \la \beta, \wt(y) \ra)\delta$ is a (not necessarily positive) real root, i.e. $w = y s_{\beta + (m + \la \beta, \wt(y) \ra)\delta}$; we call an arrow $y \xrightarrow[\beta]{m} w$ of the former type a translation edge, and one of the latter type a reflection edge. To each edge $y \xrightarrow[\beta]{m} w$, we assign an integer $d(y \xrightarrow[\beta]{m} w)$ as follows:
\begin{align}
d(y \xrightarrow[\beta]{m} w) = \begin{cases}
\frac{1}{2}\ell^\si(y,w) + m \qu & \mathrm{if}\ w = yt_{m\beta^\vee}, \\
\frac{1}{2}\bigl(\ell^\si(y,w) + 1\bigr) + m \qu & \mathrm{if}\ w = \cl(y)s_\beta t_{\wt(y) + m\beta^\vee}\ \mathrm{and}\ \cl(y) < \cl(y)s_\beta,\\
\frac{1}{2}\bigl(\ell^\si(y,w) - 1\bigr) + m \qu & \mathrm{if}\ w = \cl(y)s_\beta t_{\wt(y) + m\beta^\vee}\ \mathrm{and}\ \cl(y)s_\beta < \cl(y).
\end{cases} \nonumber
\end{align}
\item[$(2)$] A label-increasing path of length $k$ from $y$ to $w$ with respect to the reflection order $\preceq$ is a sequence $\nb = (y=y_0 \xrightarrow[\beta_{1}]{m_1} y_1 \xrightarrow[\beta_{2}]{m_2} \cdots \xrightarrow[\beta_{k}]{m_k} y_k =w)$ of elements of $\Waf$ and arrows such that $\beta_1 \prec \cdots \prec \beta_k$; we simply call this a path from $y$ to $w$. The length $k$ of a path $\nb$ is denoted by $\ell(\nb)$. Also, we denote by $P^\prec(y,w)$ the set of all paths from $y$ to $w$.
\end{enumerate}
\end{defi}

Let $\nb = (y=y_0 \xrightarrow[\beta_{1}]{m_1} y_1 \xrightarrow[\beta_{2}]{m_2} \cdots \xrightarrow[\beta_{k}]{m_k} y_k =w) \in P^\prec(y,w)$. We set
\begin{enumerate}
\item[] $e(\nb) := \{ \beta_{j} \mid j=1,\ldots,k\}$,
\item[] $r(\nb) := \{ \beta_j \mid \ \xrightarrow[\beta_j]{m_j}\ \mathrm{is\ a\ reflection\ edge} \}$,
\item[] $t(\nb) := \{ \beta_j \mid \ \xrightarrow[\beta_j]{m_j}\ \mathrm{is\ a\ translation\ edge} \}$, and
\item[] $T(\nb) := (m_{\beta})_{\beta \in \Phi_+}$, which is ordered in the (fixed) reflection order $\preceq$, where $m_{\beta}$ is $m_j$ if $\beta = \beta_j$, and is zero otherwise.
\end{enumerate}
Also, we define the degree of $\nb$ to be the integer:
\begin{align}
\deg(\nb) &= \sum_{j = 1}^{\ell(\nb)} d(y_{j-1} \xrightarrow[\beta_{j}]{m_j} y_j)- \ell(\nb).  \nonumber
\end{align}


\begin{rem}\normalfont
Since $\ell^\si(y,w)$ is even $($resp., odd$)$ if $w = yt_{m\beta^\vee}$ $($resp., $w=\cl(y)s_\beta t_{\wt(y) + m\beta^\vee}$$)$, $d(y \xrightarrow[\beta]{m} w)$ is indeed an integer.
\end{rem}



Now we are ready to state our main result of this paper.
\begin{theo}\label{main}
Let $y,w \in \Waf$. Then we have
\begin{align}
\rp_{y,w}(q) = \sum_{\nb \in P^\prec(y,w)} q^{\deg(\nb)}(q-1)^{\ell(\nb)}. \nonumber
\end{align}
\end{theo}

\begin{ex}\normalfont
In type $A_2$, let $y = e$, $w = w_0t_{\alpha_1^\vee+\alpha_2^\vee}$, $\alpha_1 < \alpha_1+\alpha_2 < \alpha_2$. Then, the elements of $P^\prec(y,w)$ are as follows:
\begin{align}
\nb_1 &= (e \xrightarrow[\alpha_1+\alpha_2]{1} w_0t_{\alpha_1^\vee+\alpha_2^\vee}), \nonumber\\
\nb_2 &= (e \xrightarrow[\alpha_1]{0} s_1 \xrightarrow[\alpha_1+\alpha_2]{1} s_2s_1t_{\alpha_1^\vee+\alpha_2^\vee} \xrightarrow[\alpha_2]{0} w_0t_{\alpha_1^\vee+\alpha_2^\vee}), \nonumber\\
\nb_3 &= (e \xrightarrow[\alpha_1]{1} t_{\alpha_1^\vee} \xrightarrow[\alpha_1+\alpha_2]{0} w_0t_{\alpha_1^\vee} \xrightarrow[\alpha_2]{1} w_0t_{\alpha_1^\vee+\alpha_2^\vee}), \nonumber\\
\nb_4 &= (e \xrightarrow[\alpha_1]{1} s_1t_{\alpha_1^\vee} \xrightarrow[\alpha_1+\alpha_2]{0} s_2s_1t_{\alpha_1^\vee} \xrightarrow[\alpha_2]{1} w_0t_{\alpha_1^\vee + \alpha_2^\vee}). \nonumber 
\end{align}
In this case, we have\\
\begin{center}
\begin{tabular}{lllll}
$r(\nb_1) = \{ \alpha_1+\alpha_2 \}$, & $t(\nb_1) = \emptyset$, & $T(\nb_1) = (0,1,0)$, & $\deg(\nb_1) = 4$, & $\ell(\nb_1) = 1$, \\
$r(\nb_2) = \Phi_+$, & $t(\nb_2) = \emptyset$, & $T(\nb_2) = (0,1,0)$, & $\deg(\nb_2) = 3$, & $\ell(\nb_2) = 3$, \\
$r(\nb_3) = \{ \alpha_1+\alpha_2 \}$, & $t(\nb_3) = \{ \alpha_1,\ \alpha_2 \}$, & $T(\nb_3) = (1,0,1)$, & $\deg(\nb_3) = 3$, & $\ell(\nb_3) = 3$, \\
$r(\nb_4) = \Phi_+$, & $t(\nb_4) = \emptyset$, & $T(\nb_4) = (1,0,1)$, & $\deg(\nb_4) = 4$, & $\ell(\nb_4) = 3$.
\end{tabular}
\end{center}
Therefore, we have
\begin{align}
\rp_{e,w_0t_{\alpha_1^\vee+\alpha_2^\vee}}(q) &= q^4(q-1) + 2q^3(q-1)^3 + q^4(q-1)^3. \nonumber
\end{align}
\end{ex}

\section{Bijections between sets of paths}
\subsection{Definitions}
We take and fix a reflection order $\preceq$ on the set $\Phi_+$ of positive roots of $\g$. For $y,w \in \Waf$, let $P^\prec_r(y,w)$ denote the subset of $P^\prec(y,w)$ that consists of all those $\nb \in P^\prec(y,w)$ having no translation edges; in the following, from an edge $u \xrightarrow[\beta]{m} v$ $(u,v \in \Waf,\ \beta \in \Phi_+)$ appearing in $\nb \in P^\prec(y,w)$, we often drop the label $m$.

\begin{defi}\normalfont
Let $u,v \in \Waf$, $s \in S_{\af}$, and let $\nb = (u = u_0 \xrightarrow[\beta_1]{} u_1 \xrightarrow[\beta_2]{} \cdots \xrightarrow[\beta_k]{} u_k = v) \in P^\prec_r(u,v)$.
\begin{enumerate}
\item[$(1)$] The $s$-descent set of $\nb$ is defined to be
\begin{align}
D_s(\nb) := \{ i \in \{ 1,2,\ldots,k \} \mid su_i = u_{i-1} \}. \nonumber
\end{align}

\item[$(2)$] The number $d_{\nb}$ $($resp., $d^\nb)$ is defined to be the minimum $($resp., maximum$)$ element of $D_s(\nb)$. If the $s$-descent set is empty, then we set $d_\nb = k+1$ and $d^\nb = 0$.
\end{enumerate}
\end{defi}

Let $P^\prec_0(u,v)$ denote the subset $\{ \nb \in P^\prec_r(u,v) \mid D_s(\nb) = \emptyset \}$ of $P^\prec_r(u,v)$, and let $P^\prec_+(u,v)$ denote the subset $\{ \nb \in P^\prec_r(u,v) \mid D_s(\nb) \neq \emptyset \}$ of $P^\prec_r(u,v)$. Now, we define a map $\Psi_L:P^\prec_+(u,v) \rightarrow P^\prec_r(su,v)$ by
\begin{align}
\Psi_L(\nb) = (su = su_0 \xrightarrow[\beta_1]{} \cdots \xrightarrow[\beta_{d_\nb - 1}]{} su_{d_\nb - 1} = u_{d_\nb} \xrightarrow[\beta_{d_\nb + 1}]{} \cdots \xrightarrow[\beta_k]{} u_k = v). \nonumber
\end{align}

Also, we define a map $\Psi_R:P^\prec_+(u,v) \rightarrow P^\prec_r(u,sv)$ by:
\begin{align}
\Psi_R(\nb) = (u = u_0 \xrightarrow[\beta_1]{} \cdots \xrightarrow[\beta_{d^\nb - 1}]{} u_{d^\nb - 1} = su_{d^\nb} \xrightarrow[\beta_{d^\nb + 1}]{} \cdots \xrightarrow[\beta_k]{} su_k = sv). \nonumber
\end{align}

\begin{rem}\normalfont\label{10.1}
We need to check that $\Psi_L(\nb) \in P^\prec_r(su,v)$ for all $\nb \in P^\prec_+(u,v)$, and that $\Psi_R(\nb) \in P^\prec_r(u,sv)$ for all $\nb \in P^\prec_+(u,v)$. For this, it suffices to show that for $u \sib v \in \Waf$, $s \in S_{\af}$, $\beta \in \Phi_+$ such that $u \xrightarrow[\beta]{m} v$ is a reflection edge, with $v \neq su$, we have $su \sib sv$. This is shown as follows. Let $\alpha$ be the simple root corresponding to the simple reflection $s$. Since $u \xrightarrow[\beta]{m} v$ is a reflection edge, we have $v = us_{\beta'}$, where $\beta' = \beta + m'\delta$ for some $m' \in \Z$; note that $v = us_{\beta'} = s_{u(\beta')}u$. Let $\gamma'$ be the unique element of $\repo \cap \{ \pm u(\beta') \}$, and write it as $\gamma' = \gamma + n\delta$ for some $n \in \Z_{\geq 0}$. Then, by Proposition $\ref{criterion}$, we have $(\cl(u))\inv(\gamma) \in \Phi_+$. Also, we have $sv = ss_{\gamma'}u = s_{s(\gamma')}su$, and $s(\gamma') = s(\gamma) + n\delta$. If $s(\gamma) \in \Phi_+$, then $s(\gamma') \in \repo$ and $(s\cl(u))\inv(s(\gamma)) = (\cl(u))\inv(\gamma) \in \Phi_+$. Hence it follows from Proposition $\ref{criterion}$ that $su \sib sv$. If $s(\gamma) \in -\Phi_+$, then $\gamma = \alpha$ and $n \neq 0$ since $\gamma' \neq \alpha$. Therefore, we obtain $s(\gamma') \in \repo$. Hence it follows from the same calculation as above that $su \sib sv$.
\end{rem}

\subsection{Some properties of paths}
By the proof of \cite[Proposition 5.2.1]{bb} (in which we order $\{ \alpha_0 \} \sqcup \Pi$ so that $\alpha_0$ is the smallest), we can show that there exists a reflection order on $\repo$ that extends the fixed reflection order $\preceq$ on $\Phi_+$; in what follows we take such a one and denote it by the same symbol. For three positive real roots $\alpha,\beta,\gamma \in \repo$, we say that $\alpha$ is between $\beta$ and $\gamma$ if either $\beta \prec \alpha \prec \gamma$ or $\gamma \prec \alpha \prec \beta$ holds.

For each $\beta \in \Phi^{\mathrm{re}}_{\af}$, we set
\begin{align}
|\beta| := \begin{cases}
\beta \qu & \mathrm{if}\ \beta \in \repo, \\
-\beta \qu & \mathrm{otherwise}.
\end{cases} \nonumber
\end{align}
Now we define two maps $\cl:\repo \rightarrow \Phi$ and $\cl_+:\Phi^{\mathrm{re}}_{\af} \rightarrow \Phi_+$ by
\begin{align}
\cl(\alpha + m\delta) &= \alpha, \nonumber\\
\cl_+(\alpha + m\delta) &= |\alpha|. \nonumber
\end{align}

\begin{lem}\label{10.2.1}
Let $u,v \in \Waf$, $s \in S_{\af}$, and $\beta \in \repo$ be such that $v = us_\beta$, and $su \sib u$. Let $\alpha$ denote the simple root corresponding to the simple reflection $s$.
\begin{enumerate}
\item[$(1)$] If $sv \sib v$, then exactly one of the following holds:
\begin{enumerate}
\item[$(\mathrm{a})$] $\cl_+(u\inv(\alpha))$ is between $\cl_+(v\inv(\alpha))$ and $\cl_+(\beta)$;
\item[$(\mathrm{b})$] $\cl_+(v\inv(\alpha))$ is between $\cl_+(u\inv(\alpha))$ and $\cl_+(\beta)$;
\item[$(\mathrm{c})$] $\cl_+(u\inv(\alpha)) = \cl_+(v\inv(\alpha)) \neq \cl_+(\beta)$.
\end{enumerate}

\item[$(2)$] If $v \sib sv$, then exactly one of the following holds:
\begin{enumerate}
\item[$(\mathrm{d})$] $\cl_+(\beta)$ is between $\cl_+(u\inv(\alpha))$ and $\cl_+(v\inv(\alpha))$;
\item[$(\mathrm{e})$] $\cl_+(u\inv(\alpha)) = \cl_+(v\inv(\alpha)) = \cl_+(\beta)$.
\end{enumerate}

\end{enumerate}
\end{lem}

\begin{proof}
$(1)$ Take $\lambda \in \Qv_+$ such that $ut_\lambda,vt_\lambda,sut_\lambda,svt_\lambda \subset \mathcal{C}^+$, and set $z := ut_\lambda, x := vt_\lambda$. Then we have
\begin{align}
sz < z, \qu sx < x, \ \text{and} \ x = zs_{\beta'}\ \text{for some}\ \beta' \in \repo, \ \text{with} \ \cl(\beta') = \cl_+(\beta). \nonumber
\end{align}
The first two conditions imply that $z\inv(\alpha), x\inv(\alpha) \in -\repo$; note that $\cl(-z\inv(\alpha)) = \cl_+(u\inv(\alpha))$ and $\cl(-x\inv(\alpha)) = \cl_+(v\inv(\alpha))$. Also, we have
\begin{align}
-z\inv(\alpha) = -s_{\beta'}x\inv(\alpha) = -x\inv(\alpha) - \la -x\inv(\alpha),(\beta')^\vee \ra \beta'. \nonumber
\end{align}
If $\la-x\inv(\alpha),(\beta')^\vee \ra < 0$, then
\begin{align}
-z\inv(\alpha) = -x\inv(\alpha) + (-\la-x\inv(\alpha),(\beta')^\vee \ra)\beta'. \nonumber
\end{align}
Applying the map $\cl$ to this equation, we obtain
\begin{align}
\cl_+(u\inv(\alpha)) = \cl_+(v\inv(\alpha)) + (-\la-x\inv(\alpha),(\beta')^\vee\ra)\cl_+(\beta). \nonumber
\end{align}
If $\cl_+(v\inv(\alpha)) = \cl_+(\beta)$, then $\cl_+(u\inv(\alpha))$ must be $\cl_+(v\inv(\alpha))$. However, this contradicts that $\la-x\inv(\alpha),(\beta')^\vee\ra \neq 0$. Hence $\cl_+(v\inv(\alpha)) \neq \cl_+(\beta)$. It follows from the definition of reflection orders that $\cl_+(u\inv(\alpha))$ is between $\cl_+(v\inv(\alpha))$ and $\cl_+(\beta)$. Thus, assertion $(a)$ holds.

If $\la -x\inv(\alpha),(\beta')^\vee \ra > 0$, then by a similar argument, assertion $(b)$ holds.

If $\la -x\inv(\alpha),(\beta')^\vee \ra = 0$, then $\cl_+(u\inv(\alpha)) = \cl_+(v\inv(\alpha))$. If, in addition, $\cl_+(u\inv(\alpha)) = \cl_+(\beta)$, then
\begin{align}
0 < \la \cl_+(u\inv(\alpha)),\cl_+(\beta) \ra = \la -x\inv(\alpha),(\beta')^\vee \ra = 0, \nonumber
\end{align}
which is a contradiction. Hence assertion $(c)$ holds.

$(2)$ We use the same notation as in the proof of $(1)$. In this case, $x\inv(\alpha) \in \repo$ and $\cl(x\inv(\alpha)) = \cl_+(v\inv(\alpha))$. Then we have
\begin{align}
-z\inv(\alpha) + x\inv(\alpha) =  - \la-x\inv(\alpha),(\beta')^\vee\ra\beta'. \nonumber
\end{align}
Since the left-hand side is a sum of positive real roots, we deduce that $a := - \la -x\inv(\alpha),(\beta')^\vee \ra > 0$. Therefore,
\begin{align}
\beta' = a\inv (-z\inv(\alpha)) + a\inv x\inv(\alpha). \nonumber
\end{align}
Applying the map $\cl$ to this equation, we obtain
\begin{align}
\cl_+(\beta) = a\inv \cl_+(u\inv(\alpha)) + a\inv \cl_+(v\inv(\alpha)). \nonumber
\end{align}
If $\cl_+(u\inv(\alpha)) \neq \cl_+(v\inv(\alpha))$, then by the definition of reflection orders, $\cl_+(\beta)$ is between $\cl_+(u\inv(\alpha))$ and $\cl_+(v\inv(\alpha))$. This proves the lemma.
\end{proof}

In order to make arguments below easier to follow, we introduce some notation.
\begin{defi}\normalfont
Let $u,v \in \Waf$, $s \in S_{\af}$, and $\beta \in \Phi_+$.
\begin{enumerate}
\item[$(1)$] We write $u \xrightarrow[\beta]{s} v$, $\xymatrix{v \\ u \ar[u]^s_{\beta}}$, or $\xymatrix{u \ar[d]_s^\beta \\ v}$ to indicate that $u \xrightarrow[\beta]{} v$ and $v = su$.

\item[$(2)$] We write $\xymatrix{v \\ u \ar@{-}[u]^s_{\beta}}$ to indicate that $v = su$ and either $u \xrightarrow[\beta]{} v$ or $v \xrightarrow[\beta]{} u$ holds.

\item[$(3)$] We write $\xymatrix{v \\ u \ar@{=}[u]}$ to indicate that $v = u$.
\end{enumerate}
\end{defi}

For example, the assumption of Lemma \ref{10.2.1} (2) is written as:
\begin{align}
\xymatrix{
u \ar@{-}[rr]_{\cl_+(\beta)} & & v \ar[dd]_s^{\cl_+(v\inv(\alpha))}\\
\\
su \ar[uu]^s_{\cl_+(u\inv(\alpha))} \ar@{-}[rr]_{\cl_+(\beta)} & & sv.
} \nonumber
\end{align}

Until the end of this subsection, we assume the following: $y \sib w \in \Waf$, $s \in S_{\af}$, $y \sib sy$, $sw \sib w$, and $u_i,v_i \in \Waf$ $(i \in \Z_{\geq0})$. For each $i \in \Z_{\geq0}$, we define $\zeta_i$ and $\xi_i$ by
\begin{align}
\zeta_i = \cl_+(u_i\inv(\alpha)), \qu \xi_i = \cl_+(v_i\inv(\alpha)), \nonumber
\end{align}
where $\alpha$ is the simple root corresponding to the simple reflection $s$.

\begin{lem}\label{10.2.3}
\begin{enumerate}
\item[$(1)$] Let $\nb = (sy = u_0 \xrightarrow[\beta_1]{} \cdots \xrightarrow[\beta_k]{} u_k = w) \in P^\prec_r(sy,w)$. Then there exists $l \in \{ 0,1,\ldots,d_\nb-1 \}$ such that $su_l \sib u_l$ and $\beta_l \prec \zeta_l \prec \beta_{l+1}$, where $\zeta_0 \prec \beta_1$ if $l = 0$, and $\beta_k \prec \zeta_k$ if $l = k$.

\item[$(2)$] Let $\nb = (y = v_0 \xrightarrow[\gamma_1]{} \cdots \xrightarrow[\gamma_k]{} v_k = sw) \in P^\prec_r(y,sw)$. Then there exists $r \in \{ d^\nb,d^\nb+1,\ldots,k \}$ such that $v_r \sib sv_r$ and $\gamma_r \prec \xi_r \prec \gamma_{r+1}$, where $\xi_0 \prec \gamma_1$ if $r = 0$, and $\gamma_k \prec \xi_k$ if $r = k$.
\end{enumerate}
\end{lem}

\begin{proof}
We prove part $(1)$; the proof of part $(2)$ is similar. Assume first that $D_s(\nb) \neq \emptyset$. What we need to do is to find $l$ for which the following diagram is commutative:
\begin{align}
\xymatrix{
sy \ar@{=}[r] & u_0 \ar[r]_{\beta_1} & u_1 \ar[r]_{\beta_2} & \cdots \ar[r]_{\beta_l} & u_l \ar[r]_{\beta_{l+1}} & \cdots \ar[r]_{\beta_{d_\nb-1}} & u_{d_\nb-1} \ar[dd]_(0.35)s^(0.35){\zeta_{d_\nb-1}} \ar[r]^s_{\beta_d} & u_{d_\nb} \\
\\
y \ar[uu]^s \ar@{=}[r] & su_0 \ar[uu]^s_{\zeta_0} \ar[r]_{\beta_1} & su_1 \ar@{-}[uu]^s_{\zeta_1} \ar[r]_{\beta_2} & \cdots \ar[r]_{\beta_l} & su_l \ar[uu]^s_{\zeta_l} \ar[r]_{\beta_{l+1}} & \cdots \ar[r]_{\beta_{d_\nb-1}} & su_{d_\nb-1}, \ar@{=}[uur]
} \nonumber
\end{align}
with $\beta_l \prec \zeta_l \prec \beta_{l+1}$. If $\zeta_0 \prec \beta_1$, then there is nothing to prove. Assume now that $\beta_1 \preceq \zeta_0$.

\begin{enumerate}
\item[]Case 1: $su_1 \sib u_1$. By Lemma \ref{10.2.1} $(1)$, we have $\beta_1 \prec \zeta_0 \prec \zeta_1$, $\beta_1 \prec \zeta_1 \prec \zeta_0$, or $\zeta_0 = \zeta_1 \neq \beta_1$. In each case, we obtain $\beta_1 \prec \zeta_1$.
\item[]Case 2: $u_1 \sib su_1$. By Lemma \ref{10.2.1} $(2)$, we have $\zeta_1 \prec \beta_1 \prec \zeta_0$, or $\beta_1 = \zeta_0 = \zeta_1$. In each case, we obtain $\zeta_1 \prec \beta_2$.
\item[]Case 2-1: $su_2 \sib u_2$. By Lemma \ref{10.2.1} $(2)$, we have $\zeta_1 \prec \beta_2 \prec \zeta_2$.
\item[]Case 2-2: $u_2 \sib su_2$. By Lemma \ref{10.2.1} $(1)$, we have $\zeta_2 \prec \zeta_1 \prec \beta_2$, $\zeta_1 \prec \zeta_2 \prec \beta_2$, or $\zeta_1 = \zeta_2 \neq \beta_2$. In each case, we obtain $\zeta_2 \prec \beta_2$, and hence $\zeta_2 \prec \beta_3$.
\end{enumerate}

If we cannot find a required $l$, then by applying the above argument repeatedly, we have either
\begin{align}
\xymatrix{
u_{d_\nb-2} \ar[r]_{\beta_{d_\nb-1}} & u_{d_\nb-1} \ar[d]_s^{\zeta_{d_\nb-1} = \beta_{d_\nb}}\\
su_{d_\nb-2} \ar[u]^s_{\zeta_{d_\nb-2}} \ar[r]_{\beta_{d_\nb-1}} & su_{d_\nb-1},
}  \nonumber
\end{align}
with $\beta_{d_\nb-1} \prec \zeta_{d_\nb-2}$, or
\begin{align}
\xymatrix{
u_{d_\nb-2} \ar[r]_{\beta_{d_\nb-1}} \ar[d]_s^{\zeta_{d_\nb-2}} & u_{d_\nb-1} \ar[d]_s^{\zeta_{d_\nb-1} = \beta_{d_\nb}}\\
su_{d_\nb-2} \ar[r]_{\beta_{d_\nb-1}} & su_{d_\nb-1},
} \nonumber
\end{align}
with $\zeta_{d_\nb-2} \prec \beta_{d_\nb-1}$. In the former case, by Lemma \ref{10.2.1} $(2)$, we see that $\zeta_{d_\nb-1} \prec \beta_{d_\nb-1} \prec \zeta_{d_\nb-2}$. Since $\zeta_{d_\nb-1} = \beta_{d_\nb}$, this contradicts the assumption that $\beta_{d_\nb-1} \prec \beta_{d_\nb}$. In the latter case, by Lemma \ref{10.2.1} $(1)$, we see that $\zeta_{d_\nb-1} \prec \zeta_{d_\nb-2} \prec \beta_{d_\nb-1}$, $\zeta_{d_\nb-2} \prec \zeta_{d_\nb-1} \prec \beta_{d_\nb-1}$, or $\zeta_{d_\nb-2} = \zeta_{d_\nb-1}$. In each case, we deduce that $\beta_{d_\nb} = \zeta_{d_\nb-1} \prec \beta_{d_\nb-1}$, which is impossible since $\beta_{d_\nb-1} \prec \beta_{d_\nb}$.

Next, we assume that $D_s(\nb) = \emptyset$. If we cannot find a required $l$ among $\{ 0,1,\ldots,k-1 \}$, then by the same argument as above, we have
\begin{align}
\xymatrix{
u_{k-1} \ar[r]_{\beta_{k}} & u_{k} \\
su_{k-1} \ar[u]^s_{\zeta_{k-1}} \ar[r]_{\beta_{k}} & su_{k} \ar[u]^s_{\zeta_{k}},
}  \nonumber
\end{align}
with $\beta_{k} \prec \zeta_{k-1}$, or
\begin{align}
\xymatrix{
u_{k-1} \ar[r]_{\beta_{k}} \ar[d]_{s}^{\zeta_{k-1}} & u_{k} \\
su_{k-1} \ar[r]_{\beta_{k}} & su_{k} \ar[u]^s_{\zeta_{k}},
} \nonumber
\end{align}
with $\zeta_{k-1} \prec \beta_{k}$. By Lemma $\ref{10.2.1}$, in each case we deduce that $\beta_{k} \prec \zeta_k$. Therefore, $l = k = d_\nb -1$ is a required one. This proves the lemma.
\end{proof}

\begin{rem}\normalfont\label{10.2.4}
Under the assumption of Lemma $\ref{10.2.3}$ $(1)$, we take an arbitrary $l$ for which the assertion of the lemma holds, and consider the path
\begin{align}
\xymatrix{
&  &  & u_l \ar[r]_{\beta_{l+1}} & \cdots \ar[r]_{\beta_{d_\nb-1}} & u_{d_\nb-1} \ar[r]_{\beta_{d_\nb}}^s  & u_{d_\nb} \ar[r]_{\beta_{d_\nb+1}} & \cdots \ar[r]_{\beta_k} & u_k = w.\\
\\
y = su_0  \ar[r]_{\beta_1} & su_1  \ar[r]_{\beta_2} & \cdots \ar[r]_{\beta_l} & su_l \ar[uu]^s_{\zeta_l}  &  & 
} \nonumber
\end{align}
Clearly, this path is an element of $\Psi_L\inv(\nb)$, where $\Psi_L:P^\prec_+(y,w) \rightarrow P^\prec(sy,w)$. Moreover, this procedure gives a bijection from the set of $l \in \{ 0,1,\ldots,d_\nb-1 \}$ such that $su_l \sib u_l$ and $\beta_l \prec \zeta_l \prec \beta_{l+1}$ to the set $\Psi_L\inv(\nb)$. Similarly, by Lemma $\ref{10.2.3}$ $(2)$, we obtain a bijection from the set of $r \in \{ d^\nb,d^\nb+1,\ldots,k \}$ such that $v_r \sib sv_r$ and $\gamma_r \prec \xi_r \prec \gamma_{r+1}$ to the set $\Psi_R\inv(\nb)$.
\end{rem}

\begin{lem}\label{10.2.5}
If we are given $(u_{l_1} \xrightarrow[\beta_{l_1 + 1}]{} u_{l_1 + 1} \xrightarrow[\beta_{l_1 + 2}]{} \cdots \xrightarrow[\beta_{l_2}]{} u_{l_2}) \in P_0^<(u_{l_1},u_{l_2})$ such that $su_{l_1} \sib u_{l_1}$, $su_{l_2} \sib u_{l_2}$, $\zeta_{l_1} \prec \beta_{l_1 + 1}$, and $\beta_{l_2} \prec \zeta_{l_2}$, then we can find $d \in \{ l_1 + 1,\ldots,l_2-1\}$ such that $u_{d} \sib su_{d}$ and $\beta_{d} \prec \zeta_{d} \prec \beta_{d+1}$.
\end{lem}

\begin{proof}
We want to find $d$ for which the following diagram is commutative:
\begin{align}
\xymatrix{
u_{l_1} \ar[r]_{\beta_{l_1+1}} & \cdots \ar[r]_{\beta_{d}} & u_{d} \ar[r]_{\beta_{d+1}} \ar[dd]_{s}^{\zeta_{d}} & \cdots \ar[r]_{\beta_{l_2}} & u_{l_2} \\
\\
su_{l_1} \ar[uu]^s_{\zeta_{l_1}} \ar[r]_{\beta_{l_1+1}} & \cdots \ar[r]_{\beta_{d}} & su_{d} \ar[r]_{\beta_{d+1}} & \cdots \ar[r]_{\beta_{l_2}} & su_{l_2}, \ar[uu]_{\zeta_{l_2}}^s
} \nonumber
\end{align}
with $\beta_{d} \prec \zeta_{d} \prec \beta_{d+1}$. First, we prove that there exists $d'$ such that $u_{d'} \sib su_{d'}$. If this does not hold, then we have (by Lemma \ref{10.2.1} $(1)$) $\zeta_{l_1 + 1} \prec \beta_{l_1 + 1} (\prec \beta_{l_1 +2})$, and inductively
\begin{align}
\zeta_{l_1 + 2} <\ &\beta_{l_1 + 2} \prec \beta_{l_1 +3}, \nonumber \\
\zeta_{l_1 + 3} <\ &\beta_{l_1 + 3} \prec \beta_{l_1 +4}, \nonumber \\
&\hspace{0.2cm}\vdots \nonumber\\
\zeta_{l_2} <\ &\beta_{l_2}. \nonumber
\end{align}
However, the last inequality contradicts the assumption. Now, we proceed by induction on $l_2 - l_1$. If $l_2 - l_1 = 2$, then we obtain the commutative diagram:
\begin{align}
\xymatrix{
u_{l_1} \ar[r]_{\beta_{l_1+1}} & u_{l_1+1} \ar[dd]_s^{\zeta_{l_1+1}} \ar[r]_{\beta_{l_1+2}} & u_{l_1+2} \\
\\
su_{l_1} \ar[uu]^s_{\zeta_{l_1}} \ar[r]_{\beta_{l_1+1}} & su_{l_i+1} \ar[r]_{\beta_{l_1+2}} & su_{l_1+2}. \ar[uu]_{\zeta_{l_1+2}}^s
} \nonumber
\end{align}
Note that by Lemma \ref{10.2.1} $(2)$, we have $\zeta_{l_1} \prec \beta_{l_1 + 1} \prec \zeta_{l_1+1} \prec \beta_{l_1+2} \prec \zeta_{l_1+2}$. Therefore, $d := l_1+1$ is a required one. When $l_2 - l_1 > 2$, take the minimum index $d'$ such that $u_{d'} \sib su_{d'}$; it is easy to show (by using Lemma \ref{10.2.1}) that $\beta_{d'} \prec \zeta_{d'}$. If $\zeta_{d'} \prec \beta_{d'+1}$, then $d := d'$ is a required one. If not, take the minimum index $l' \in \{d'+1,\ldots,l_2 \}$ such that $su_{l'} \sib u_{l'}$. Then, by using Lemma \ref{10.2.1} again, we deduce that either $\zeta_{l'} \prec \beta_{l'}$, or there exists $d'' \in \{ d'+1,\ldots,l'-1 \}$ such that $u_{d''} \sib su_{d''}$ and $\beta_{d''} \prec \zeta_{d''} \prec \beta_{d''+1}$. Here, in the former case, we obtain $l' \neq l_2$, and hence we can complete the proof by our induction hypothesis.
\end{proof}


\begin{prop}\label{10.2.7}
Let $\nb \in P^\prec_+(y,w)$. Then the following hold.
\begin{enumerate}
\item[$(1)$] $\left| \Psi_L\inv(\nb) \right| = \left|\{ \nb' \in \Psi_L\inv(\Psi_L(\nb)) \mid d_{\nb'} < d_{\nb} \} \right|$.
\item[$(2)$] $\left|\Psi_R\inv(\nb) \right| = \left|\{ \nb' \in \Psi_R\inv(\Psi_R(\nb)) \mid d^{\nb'} > d^{\nb} \} \right|$.
\end{enumerate}
\end{prop}
\begin{proof}
We prove part $(1)$; the proof of part $(2)$ is similar. Let us write
\begin{align}
\Psi_L(\nb) = (sy = u_0 \xrightarrow[\beta_1]{} u_1 \xrightarrow[\beta_2]{} \cdots \xrightarrow[\beta_k]{} u_k = w), \nonumber
\end{align}
and $\{ \nb' \in \Psi_L\inv(\Psi_L(\nb)) \mid d_{\nb'} \leq d_{\nb} \} = \{ \nb_1,\ldots,\nb_p \}$, with $d_{\nb_1} < \cdots < d_{\nb_p}$. Observe that if $\Psi_L(\nb') = \Psi_L(\nb'')$ and $d_{\nb'} = d_{\nb''}$, then $\nb' = \nb''$. It follows from Remark $\ref{10.2.4}$ that $\nb_p = \nb$. Also, from Lemma \ref{10.2.3} (1) and Remark \ref{10.2.4}, we have
\begin{align}
su_{d_{\nb_i}-1} \sib u_{d_{\nb_i}-1}\ and\ \beta_{d_{\nb_i}-1} \prec \zeta_{d_{\nb_i}-1} \prec \beta_{d_{\nb_i}} \nonumber
\end{align}
for all $i = 1,\ldots,p$. Let $i \geq 2$. Then the path
\begin{align}
(u_{d_{\nb_{i-1}}-1} \xrightarrow[\beta_{d_{\nb_{i-1}}}]{} \cdots \xrightarrow[\beta_{d_{\nb_{i}}-1}]{} u_{d_{\nb_{i}}-1}) \nonumber
\end{align}
satisfies the assumption of Lemma \ref{10.2.5}. Therefore, there exists $d_i \in \{ d_{\nb_{i-1}},\ldots,d_{\nb_i}-2 \}$ for which the following diagram is commutative:
\begin{align}
\xymatrix{
\cdots \ar[r]_{\beta_{d_{\nb_{i-1}}-1}} &u_{d_{\nb_{i-1}}-1} \ar[r]_{{\beta_{d_{\nb_{i-1}}}}} & \cdots \ar[r]_{\beta_{d_i}} & u_{d_i} \ar[dd]_s^{\zeta_{d_i}} \ar[r]_{\beta_{d_i+1}} & \cdots \ar[r]_{\beta_{d_{\nb_i}-1}} & u_{d_{\nb_i}-1} \ar[r]_{\beta_{d_{\nb_i}}} & \cdots \ar[r]_{\beta_{d_{\nb_p}-1}} & u_{d_{\nb_p}-1} \ar[r]_{\beta_{d_{\nb_p}}} & \cdots \\
\\
\cdots \ar[r]_{\beta_{d_{\nb_{i-1}}-1}} &su_{d_{\nb_{i-1}}-1} \ar[uu]^s_{\zeta_{d_{\nb_{i-1}}-1}} \ar[r]_{{\beta_{d_{\nb_{i-1}}}}} & \cdots \ar[r]_{\beta_{d_i}} & su_{d_i} \ar[r]_{\beta_{d_i+1}} & \cdots \ar[r]_{\beta_{d_{\nb_i}-1}} & su_{d_{\nb_i}-1} \ar[uu]^s_{\zeta_{d_{\nb_i}-1}} \ar[r]_{\beta_{d_{\nb_i}}} & \cdots \ar[r]_{\beta_{d_{\nb_p}-1}} & su_{d_{\nb_p}-1} \ar[uu]^s_{\zeta_{d_{\nb_p}-1}} \ar[r]_{\beta_{d_{\nb_p}}} & \cdots,
} \nonumber
\end{align}
with $\beta_{d_i} \prec \zeta_{d_i} \prec \beta_{d_i + 1}$. If we set
\begin{align}
\nb'_i &:= (sy = u_0 \xrightarrow[\beta_1]{} \cdots \xrightarrow[\beta_{d_i}]{} u_{d_i} \xrightarrow[\zeta_{d_i}]{s} su_{d_i} \xrightarrow[\beta_{d_i+1}]{} \cdots \nonumber\\
&\hspace{1cm}\xrightarrow[\beta_{d_{\nb_p}-1}]{} su_{d_{\nb_p}-1} \xrightarrow[\zeta_{d_{\nb_p}-1}]{s} u_{d_{\nb_p}-1} \xrightarrow[\beta_{d_{\nb_p}}]{} \cdots \xrightarrow[\beta_k]{} u_k = w) \in P^\prec(sy,w), \nonumber
\end{align}
then we have $\nb'_i \in \Psi_L\inv(\nb_p) = \Psi_L\inv(\nb)$, with $d_{\nb'_i} = d_i + 1$. If there exists another $d'_i \in \{ d_{\nb_{i-1}},\ldots,d_{\nb_i}-2 \}$ such that $u_{d'_i} \sib su_{d'_i}$ and $\beta_{d'_i} \prec \zeta_{d'_i} \prec \beta_{d'_i+1}$, then we may assume that $d'_i < d_i$. By applying Lemma \ref{10.2.5} to the path $(su_{d'_i} \xrightarrow[\beta_{d'_i+1}]{} \cdots \xrightarrow[\beta_{d_i}]{} su_{d_i})$, we obtain $\nb' \in \Psi_L\inv(\Psi_L(\nb))$ such that $d_{\nb_{i-1}} < d'_i + 1 < d_{\nb'} < d_i+1 < d_{\nb_i} \leq d_{\nb_p}$. However, this is impossible because $\nb' \in \{ \nb' \in \Psi_L\inv(\Psi_L(\nb)) \mid d_{\nb'} \leq d_{\nb} \} = \{ \nb_1,\ldots,\nb_p \}$ and $d_{\nb_1} < \cdots < d_{\nb_p}$. This proves the proposition.
\end{proof}

\subsection{Bijections}
Recall that we have assumed that $y \sib w \in \Waf$, $s \in S_{\af}$, $y \sib sy$, $sw \sib w$, and $u_i,v_i \in \Waf$ $(i \in \Z_{\geq0})$. Also, for each $i \in \Z_{\geq0}$, $\zeta_i$ and $\xi_i$ are defined by
\begin{align}
\zeta_i = \cl_+(u_i\inv(\alpha)), \qu \xi_i = \cl_+(v_i\inv(\alpha)), \nonumber
\end{align}
where $\alpha$ is the simple root corresponding to the simple reflection $s$. In the following, we define a map $\varphi : P^\prec_r(sy,sw) \sqcup P^\prec_r(sy,w) \rightarrow P^\prec_r(y,w)$. If $\nb = (sy = u_0 \xrightarrow[\beta_1]{} u_1 \xrightarrow[\beta_2]{} \cdots \xrightarrow[\beta_k]{} u_k = sw) \in P^\prec_0(sy,sw)$, then we define
\begin{align}
\varphi(\nb) := (y = su_0 \xrightarrow[\beta_1]{} su_1 \xrightarrow[\beta_2]{} \cdots \xrightarrow[\beta_k]{} su_k = w); \nonumber
\end{align}
note that $\ell(\varphi(\nb)) = \ell(\nb)$.

Let $\nb = (sy = u_0 \xrightarrow[\beta_1]{} u_1 \xrightarrow[\beta_2]{} \cdots \xrightarrow[\beta_k]{} u_k = w) \in P^\prec_r(sy,w)$. By Remark \ref{10.2.4}, the set $\Psi_L\inv(\nb)$ is not empty. We define $\varphi(\nb)$ to be the unique element of $\Psi_L\inv(\nb)$ such that $d_{\varphi(\nb)} < d_{\nb'}$ for all $\nb' \in \Psi_L\inv(\nb) \setminus \{ \varphi(\nb) \}$; this $\varphi(\nb)$ is realized by taking the smallest $l$ for which the assertion of Lemma \ref{10.2.3} $(1)$ holds. Note that in this case, we have $\ell(\varphi(\nb)) = \ell(\nb) + 1$.

Let $\nb = (sy = v_0 \xrightarrow[\gamma_1]{} v_1 \xrightarrow[\gamma_2]{} \cdots \xrightarrow[\gamma_k]{} v_k = sw) \in P^\prec_{+}(sy,sw)$. Let us write $D_s(\nb) = \{ d_1,\ldots,d_m \}$, with $d_1 < \cdots < d_m$. For each $i = 1,\ldots,m$, we have
\begin{align}
\xymatrix{
v_{d_i-1} \ar@{=}[r] \ar[dd]_s^{\gamma_{d_i}} & sv_{d_i} \ar[r]_{\gamma_{d_i+1}} \ar[dd]_s^{\gamma_{d_i}} & sv_{d_i+1} \ar[r]_{\gamma_{d_i+2}} & \cdots \ar[r]_{\gamma_{d_{i+1}-1}} & sv_{d_{i+1}-1} \ar@{=}[r] & v_{d_{i+1}}\\
\\
sv_{d_i-1} \ar@{=}[r] & v_{d_i} \ar[r]_{\gamma_{d_i+1}} & v_{d_i+1} \ar[r]_{\gamma_{d_i+2}} \ar@{-}[uu]^s_{\xi_{d_i+1}} & \cdots \ar[r]_{\gamma_{d_{i+1}-1}} & v_{d_{i+1}-1} \ar[uu]^s_{\gamma_{d_{i+1}}} \ar@{=}[r] & sv_{d_{i+1}}, \ar[uu]^s_{\gamma_{d_{i+1}}}
} \nonumber
\end{align}
where $d_{m+1} = k+1$ and $v_{k+1} = w$. By applying Lemma \ref{10.2.3} $(2)$ to the path $(sv_{d_i} \xrightarrow[\gamma_{d_i}]{s} v_{d_i} \xrightarrow[\gamma_{d_i+1}]{} \cdots \xrightarrow[\gamma_{d_{i+1}-1}]{} v_{d_{i+1}-1})$, we can take the smallest $r_i \in \{ d_i+1,d_i+2,\ldots,d_{i+1}-1 \}$ such that $v_{r_i} \sib sv_{r_i}$ and $\gamma_{r_i} \prec \xi_{r_i} \prec \gamma_{r_i+1}$. Then, we have
\begin{align}
\nb_i := (v_{d_i} \xrightarrow[\gamma_{d_{i}+1}]{} \cdots \xrightarrow[\gamma_{r_i}]{} v_{r_i} \xrightarrow[\xi_{r_i}]{s} sv_{r_i} \xrightarrow[\gamma_{r_i+1}]{} \cdots \xrightarrow[\gamma_{d_{i+1}-1}]{} sv_{d_{i+1}-1} = v_{d_{i+1}}) \in P^\prec_r(v_{d_i},v_{d_{i+1}}). \nonumber
\end{align}
Also, we set $\nb_0 := (y = sv_0 \xrightarrow[\gamma_1]{} sv_1 \xrightarrow[\gamma_2]{} \cdots \xrightarrow[\gamma_{d_1-1}]{} sv_{d_1-1} = v_{d_1}) \in P^\prec_0(y,v_{d_1})$. Now we define $\varphi(\nb)$ to be the concatenation of the paths $\nb_0,\nb_1,\ldots,\nb_m$. Namely, we define $\varphi(\nb) = (y = u_0 \xrightarrow[\beta_1]{} u_1 \xrightarrow[\beta_2]{} \cdots \xrightarrow[\beta_k]{} u_k = w)$, as follows:
\begin{align}
u_j &= 
\begin{cases}
sv_j \qu & (j = 0,1,\ldots,d_1-1), \\
v_{j+1} \qu & (j = d_i,d_i+1,\ldots,r_i-1,\ \mathrm{and}\ i = 1,\ldots,m), \\
sv_j \qu & (j = r_i,r_i+1,\ldots,d_{i+1}-1,\ \mathrm{and}\ i= 1,\ldots,m), 
\end{cases} \label{eq:3}\\
\beta_j &=  
\begin{cases}
\gamma_j \qu & (j = 1,\ldots,d_1-1), \\
\gamma_{j+1} \qu & (j = d_i,d_i+1,\ldots,r_i-1,\ \mathrm{and}\ i = 1,\ldots,m), \\
\gamma_j \qu & (j = r_i+1,\ldots,d_{i+1}-1,\ \mathrm{and}\ i= 1,\ldots,m), \\
\xi_{r_i} \qu & (j = r_i,\ \mathrm{and}\ i=1,\ldots,m).
\end{cases} \label{eq:4}
\end{align}
Note that in this case, we have $\ell(\varphi(\nb)) = \ell(\nb)$.

\begin{prop}\label{10.3.1}
The map $\varphi : P^\prec_r(sy,sw) \sqcup P^\prec_r(sy,w) \rightarrow P^\prec_r(y,w)$ is a bijection.
\end{prop}

\begin{proof}
We prove the assertion by constructing a map $\varphi':P^\prec_r(y,w) \rightarrow P^\prec_r(sy,sw) \sqcup P^\prec_r(sy,w)$, which turns out to be the inverse of $\varphi$. Let $\nb' = (y = u_0 \xrightarrow[\beta_1]{} \cdots \xrightarrow[\beta_k]{} u_k = w) \in P^\prec_r(y,w)$. If $\nb' \in P^\prec_0(y,w)$, then $\nb' = \varphi(\nb)$, where
\begin{align}
\nb := (sy = su_0 \xrightarrow[\beta_1]{} su_1 \xrightarrow[\beta_2]{} \cdots \xrightarrow[\beta_k]{} su_k = sw) \in P^\prec_0(sy,sw). \nonumber
\end{align}
Hence we define $\varphi'(\nb') := \nb$.

Assume that $\nb' \in P^\prec_{+}(y,w)$. If the set $\Psi_L\inv(\nb')$ is empty, then by Proposition \ref{10.2.7}, $d_{\nb'} \leq d_{\nb''}$ for all $\nb'' \in \Psi_L\inv(\Psi_L(\nb'))$. Therefore, by the definition of $\varphi$, we have
\begin{align}
\nb' = \varphi(\Psi_L(\nb')). \nonumber
\end{align}
Hence we define $\varphi'(\nb') := \Psi_L(\nb')$.

Next, assume that $\nb' \in P^\prec_{+}(y,w)$ and $\Psi_L\inv(\nb')$ is not empty. Let us write $D_s(\nb') = \{ d'_1,\ldots,d'_m \}$, with $d'_1 < \cdots < d'_m$. Since $\Psi_L\inv(\nb')$ is not empty, we can take the maximum index $l_1 \in \{ 1,\ldots,d'_1-2 \}$ such that $su_{l_1} < u_{l_1}$ and $\beta_{l_1} \prec \zeta_{l_1} \prec \beta_{l_1+1}$ (see Remark \ref{10.2.4}). For each $i = 1,\ldots,m-1$, we have the following commutative diagram:
\begin{align}
\xymatrix{
u_{d'_i-1} \ar@{=}[r] \ar[dd]_s^{\beta_{d'_i}} & su_{d'_i} \ar[r]_{\beta_{d'_i+1}} \ar[dd]_s^{\beta_{d'_i}} & su_{d'_i+1} \ar[r]_{\beta_{d'_i+2}} & \cdots \ar[r]_{\beta_{d'_{i+1}-1}} & su_{d'_{i+1}-1} \ar@{=}[r] & u_{d'_{i+1}}\\
\\
su_{d'_i-1} \ar@{=}[r] & u_{d'_i} \ar[r]_{\beta_{d'_i+1}} & u_{d'_i+1} \ar[r]_{\beta_{d'_i+2}} \ar@{-}[uu]^s_{\zeta_{d'_i+1}} & \cdots \ar[r]_{\beta_{d'_{i+1}-1}} & u_{d'_{i+1}-1} \ar[uu]^s_{\beta_{d'_{i+1}}} \ar@{=}[r] & su_{d'_{i+1}}. \ar[uu]^s_{\beta_{d'_{i+1}}}
} \nonumber
\end{align}
By applying Lemma \ref{10.2.3} $(1)$ to the path $(su_{d'_i-1} = u_{d'_i} \xrightarrow[\beta_{d'_i+1}]{} \cdots \xrightarrow[\beta_{d'_{i+1}-1}]{} u_{d'_{i+1}-1} \xrightarrow[\beta_{d'_{i+1}}]{s} su_{d'_{i+1}-1} = u_{d'_{i+1}})$, we can take the maximum index $l_{i+1} \in \{ d'_i,d'_i+1,\ldots,d'_{i+1}-2 \}$ such that $su_{l_{i+1}} < u_{l_{i+1}}$ and $\beta_{l_{i+1}} \prec \zeta_{l_{i+1}} \prec \beta_{l_{i+1}+1}$. Then, we have
\begin{align}
\nb'_i := (u_{d'_i-1} = su_{d'_i} \xrightarrow[\beta_{d'_{i}+1}]{} \cdots \xrightarrow[\beta_{d'_{l_i}}]{} su_{l_{i+1}} \xrightarrow[\zeta_{l_{i+1}}]{s} u_{l_{i+1}} \xrightarrow[\beta_{l_i+1}]{} \cdots \xrightarrow[\beta_{d'_{i+1}-1}]{} u_{d'_{i+1}-1} ) \in P^\prec_r(u_{d'_i-1},u_{d'_{i+1}-1}). \nonumber
\end{align}
Also, we set $\nb'_0 = (sy = su_0 \xrightarrow[\beta_1]{} \cdots \xrightarrow[\beta_{l_1}]{} su_{l_1} \xrightarrow[\zeta_{l_1}]{s} u_{l_1} \xrightarrow[\beta_{l_1+1}]{} \cdots \xrightarrow[\beta_{d'_1-1}]{} u_{d'_1-1}) \in P^\prec_r(sy,u_{d'_1-1})$ and $\nb'_m := (u_{d'_m-1} = su_{d'_m} \xrightarrow[\beta_{d'_m+1}]{} \cdots \xrightarrow[\beta_k]{} su_k = sw) \in P^\prec_0(u_{d'_m-1},sw)$. We define $\varphi'(\nb')$ to be the concatenation $\nb \in P^\prec_+(sy,sw)$ of the paths $\nb'_0,\nb'_1,\ldots,\nb'_m$. Namely, we define $\nb = (sy = v_0 \xrightarrow[\gamma_1]{} v_1 \xrightarrow[\gamma_2]{} \cdots \xrightarrow[\gamma_k]{} v_k = sw)$ as follows:
\begin{align}
v_j &= \begin{cases}
su_j \qu & (j = d'_i,d'_i+1,\ldots,l_{i+1},\ \mathrm{and}\ i=0,1,\ldots,m-1),\\
u_{j-1} \qu & (j = l_{i+1}+1,l_{i+1}+2,\ldots,d'_{i+1}-1,\ \mathrm{and}\ i = 0,1,\ldots,m-1), \\
su_j \qu & (j = d'_m,d'_m+1,\ldots,k),
\end{cases} \label{eq:5} \\
\gamma_j &= \begin{cases}
\beta_j \qu & (j = d'_i+1,\ldots,l_{i+1},\ \mathrm{and}\ i=0,1,\ldots,m-1),\\
\beta_{j-1} \qu & (j = l_{i+1}+2,\ldots,d'_{i+1},\ \mathrm{and}\ i = 0,1,\ldots,m-1), \\
\beta_j \qu & (j = d'_m,d'_m+1,\ldots,k), \\
\zeta_{l_{i+1}} \qu & (j = l_{i+1}+1, \ \mathrm{and}\ i = 0,1,\ldots,m-1),
\end{cases} \label{eq:6}
\end{align}
where $d'_0 := 0$.

We show that $\varphi(\nb) = \nb'$. By the definition of $\nb$, we have $D_s(\nb) = \{ l_1+1,l_2+1,\ldots,l_m+1 \}$, with $l_1 \leq d'_1 -2 < d'_1 \leq l_2 \leq d'_2 -2 < d'_2 \leq \cdots \leq l_m \leq d'_m -2$. We claim that $d'_i$ coincides with the minimum index $r_i \in \{ l_i+2,l_i+3,\ldots,l_{i+1} \}$ such that $v_{r_i} \sib sv_{r_i}$ and $\gamma_{r_i} \prec \xi_{r_i} \prec \gamma_{r_i+1}$ for all $i = 1,\ldots,m$ (here we set $l_{m+1} = k+1$). If this does not hold, then we obtain the following commutative diagram:
\begin{align}
\xymatrix{
v_{l_i} \ar@{=}[r] \ar[dd]_s & sv_{l_i+1} \ar[r] \ar[dd]_s & \cdots \ar[r] & sv_{r_i} \ar[r] & \cdots \ar[r] & sv_{d'_i} \ar[r] & \cdots \ar[r] & sv_{l_{i+1}} \ar@{=}[r] & v_{l_{i+1} + 1}\\
\\
sv_{l_i} \ar@{=}[r] & v_{l_i+1} \ar[r] & \cdots \ar[r] & v_{r_i} \ar[r] \ar[uu]^s_{\xi_{r_i}} & \cdots \ar[r] & v_{d'_i} \ar[r] \ar[uu]^s_{\beta_{d'_{i}}} & \cdots \ar[r] & v_{l_{i+1}} \ar@{=}[r] \ar[uu]^s & sv_{l_{i+1}+1}. \ar[uu]^s
} \nonumber
\end{align}
By Lemma \ref{10.2.5}, there exists $d \in \{ r_i+1,\ldots,d'_i-1 \}$ such that $su_{d-1} = sv_d \sib v_d = u_{d-1}$ and $\beta_{d-1} = \gamma_d \prec \xi_d (= \zeta_{d-1}) \prec \gamma_{d+1} \beta_{d}$. However, this contradicts the maximality of $l_i$. If we write $\varphi(\nb) = (y = u'_0 \xrightarrow[]{} u'_1 \xrightarrow[]{} \cdots \xrightarrow[]{} u'_k = w)$, then it is easy to see that $\varphi(\nb) = \nb'$ by using equations \eqref{eq:3}, \eqref{eq:4} (with $d_i = l_i + 1$, $r_i = d'_i$ for $i = 1,\ldots,m$), \eqref{eq:5}, and \eqref{eq:6}.

It remains to show that $\varphi'(\varphi(\nb)) = \nb$ for all $\nb \in P^\prec_r(sy,sw) \sqcup P^\prec_r(sy,w)$. It is easy when $\nb \in P^\prec_0(sy,sw) \sqcup P^\prec_r(sy,w)$. When $\nb \in P^\prec_+(sy,sw)$, we deduce from the definition of $\varphi(\nb)$ (in particular, that of $\nb_0$) that $\Psi_L\inv(\varphi(\nb))$ is not empty. Based on this fact, we can show that $\varphi'(\varphi(\nb)) = \nb$ in a way similar to the above. This completes the proof of the proposition.
\end{proof}

Now, we define a map $\psi:P^\prec_r(y,sw) \rightarrow P^\prec_r(sy,w)$ as follows. Let $\nb = (y = v_0 \xrightarrow[\gamma_1]{} \cdots \xrightarrow[\gamma_k]{} v_k = sw) \in P^\prec_r(y,sw)$. When $\nb \in P^\prec_0(y,sw)$, we define
\begin{align}
\psi(\nb) := (sy = sv_0 \xrightarrow[\gamma_1]{} sv_1 \xrightarrow[\gamma_2]{} \cdots \xrightarrow[\gamma_k]{} sv_k = w) \in P^\prec_0(sy,w). \nonumber
\end{align}
Assume that $\nb \in P^\prec_+(y,sw)$. We write $D_s(\nb) = \{ d_1,\ldots,d_m \}$, with $d_1 < \cdots < d_m$. As in the definition of $\varphi$, for each $i = 1,\ldots,m$, we can take the minimum index $r_i \in \{ d_i+1,\ldots,d_{i+1}-1 \}$ such that $v_{r_i} \sib sv_{r_i}$ and $\gamma_{r_i} \prec \xi_{r_i} \prec \gamma_{r_i+1}$; hence we can define $\nb_i \in P^\prec_r(v_{d_i},v_{d_{i+1}})$ by
\begin{align}
\nb_i := (v_{d_i} \xrightarrow[\gamma_{d_{i}+1}]{} \cdots \xrightarrow[\gamma_{r_i}]{} v_{r_i} \xrightarrow[\xi_{r_i}]{s} sv_{r_i} \xrightarrow[\gamma_{r_i+1}]{} \cdots \xrightarrow[\gamma_{d_{i+1}-1}]{} sv_{d_{i+1}-1} = v_{d_{i+1}}) \in P^\prec_r(v_{d_i},v_{d_{i+1}}), \nonumber
\end{align}
where $d_{m+1} = k+1$ and $v_{k+1} = w$. Also, we set $\nb_0 := (sy = sv_0 \xrightarrow[\gamma_1]{} sv_1 \xrightarrow[\gamma_2]{} \cdots \xrightarrow[\gamma_{d_1-1}]{} sv_{d_1-1} = v_{d_1}) \in P^\prec_0(sy,v_{d_1})$. We define $\psi(\nb)$ to be the concatenation of the paths $\nb_0,\nb_1,\ldots,\nb_m$.

In a way similar to the definition of $\varphi'$ above, we define a map $\psi':P^\prec_r(sy,w) \rightarrow P^\prec_r(y,sw)$. Let $\nb' = (sy = u_0 \xrightarrow[\beta_1]{} \cdots \xrightarrow[\beta_k]{} u_k = w) \in P^\prec_r(sy,w)$. If $\nb' \in P^\prec_0(sy,w)$, then $\nb' = \psi(\nb)$, where
\begin{align}
\nb := (y = su_0 \xrightarrow[\beta_1]{} \cdots \xrightarrow[\beta_k]{} su_k = sw) \in P^\prec_0(y,sw).
\label{eq:7}\end{align}
Hence we define $\psi'(\nb') := \nb$.

Assume that $\nb' \in P^\prec_+(sy,w)$. By Lemma \ref{10.2.3} and Remark \ref{10.2.4}, $\Psi_L\inv(\nb')$ is not empty. Let us write $D_s(\nb') = \{ d'_1,\ldots,d'_m \}$, with $0 = d'_0 < d'_1 < \cdots d'_m$. For each $i = 0,1,\ldots,m-1$, we have the following commutative diagram:
\begin{align}
\xymatrix{
u_{d'_i-1} \ar@{=}[r] \ar[dd]_s^{\beta_{d'_i}} & su_{d'_i} \ar[r]_{\beta_{d'_i+1}} \ar[dd]_s^{\beta_{d'_i}} & su_{d'_i+1} \ar[r]_{\gamma_{d'_i+2}} & \cdots \ar[r]_{\beta_{d'_{i+1}-1}} & su_{d'_{i+1}-1} \ar@{=}[r] & u_{d'_{i+1}}\\
\\
su_{d'_i-1} \ar@{=}[r] & u_{d'_i} \ar[r]_{\gamma_{d'_i+1}} & u_{d'_i+1} \ar[r]_{\beta_{d'_i+2}} \ar@{-}[uu]^s_{\zeta_{d'_i+1}} & \cdots \ar[r]_{\beta_{d'_{i+1}-1}} & u_{d'_{i+1}-1} \ar[uu]^s_{\beta_{d'_{i+1}}} \ar@{=}[r] & su_{d'_{i+1}}, \ar[uu]^s_{\beta_{d'_{i+1}}}
} \nonumber
\end{align}
where $u_{-1} = y$. By applying Lemma \ref{10.2.3} $(1)$ to the path $(su_{d'_i-1} = u_{d'_i} \xrightarrow[\beta_{d'_i+1}]{} \cdots \xrightarrow[\beta_{d'_{i+1}-1}]{} u_{d'_{i+1}-1} \xrightarrow[\beta_{d'_{i+1}}]{s} su_{d'_{i+1}-1} = u_{d'_{i+1}})$, we can take the maximum index $l_{i+1} \in \{ d'_i,d'_i+1,\ldots,d'_{i+1}-2 \}$ such that $su_{l_{i+1}} < u_{l_{i+1}}$ and $\gamma_{l_{i+1}} \prec \gamma_{l_{i+1}} \prec \gamma_{l_{i+1}+1}$. Then, we have
\begin{align}
\nb'_i := (u_{d'_i-1} = su_{d'_i} \xrightarrow[\beta_{d'_{i}+1}]{} \cdots \xrightarrow[\gamma_{d'_{l_i}}]{} su_{l_{i+1}} \xrightarrow[\zeta_{l_{i+1}}]{s} u_{l_{i+1}} \xrightarrow[\beta_{l_i+1}]{} \cdots \xrightarrow[\beta_{d'_{i+1}-1}]{} u_{d'_{i+1}-1} ) \in P^\prec_r(u_{d'_i-1},u_{d'_{i+1}-1}). \nonumber
\end{align}
Also, we set $\nb'_m := (u_{d'_m-1} = su_{d'_m} \xrightarrow[\beta_{d'_m+1}]{} \cdots \xrightarrow[\beta_k]{} su_k = sw) \in P^\prec_0(u_{d'_m-1},sw)$. We define $\varphi'(\nb')$ to be the concatenation $\nb \in P^\prec_+(y,sw)$ of the paths $\nb'_0,\nb'_1,\ldots,\nb'_m$.

The proof of the following proposition is similar to that of Proposition \ref{10.3.1}.

\begin{prop}
The map $\psi:P^\prec_r(y,sw) \rightarrow P^\prec_r(sy,w)$ is a bijection, with $\psi'$ its inverse.
\end{prop}

We now summarize the results obtained in this section:

\begin{cor}\label{10.3.3}
Let $y,w \in \Waf$, $s \in S_{\af}$. Assume that $sw \sib w$. Then, there exist bijections $\varphi:P^\prec_r(sy,sw) \sqcup P^\prec_r(sy,w) \rightarrow P^\prec_r(y,w)$ and $\psi:P^\prec_r(y,sw) \rightarrow P^\prec_r(sy,w)$. Moreover, they have the following properties:
\begin{align}
\ell(\varphi(\nb)) &= \begin{cases}
\ell(\nb) \qu & (\nb \in P^\prec_r(sy,sw)), \\
\ell(\nb) + 1 \qu & (\nb \in P^\prec_r(sy,w)),
\end{cases} \nonumber\\
\ell(\psi(\nb)) &= \ell(\nb). \nonumber
\end{align}
\end{cor}

\section{Proof of Theorem \ref{main}}
\subsection{Recursion relations}
Let $y,w \in \Waf$, and let $\preceq$ be a (fixed) reflection order on $\Phi_+$. Recall that $P^\prec_r(y,w)$ is defined to be the set of all paths from $y$ to $w$ without translation edges. Also, we define $P^\prec_t(y,w)$ to be the set of all paths from $y$ to $w$ without reflection edges.

\begin{defi}\normalfont
For $y,w \in \Waf$, we define
\begin{align}
r^\prec_{y,w}(q) &:= \sum_{\nb \in P^\prec_r(y,w)} q^{\frac{1}{2}\bigl(\ell^\si(y,w)-\ell(\nb)\bigr)}(q-1)^{\ell(\nb)}, \nonumber\\
t^\prec_{y,w}(q) &:= \sum_{\nb \in P^\prec_t(y,w)} q^{\deg(\nb)}(q-1)^{\ell(\nb)}. \nonumber
\end{align}
\end{defi}

\begin{prop}\label{11.1.2}
Let $y,w \in \Waf$, and $s \in S_{\af}$. Assume that $sw\sib w$.
\begin{enumerate}
\item[$(1)$] It holds that \begin{align}
r^\prec_{y,w}(q) &= \begin{cases}
r^\prec_{sy,sw}(q) \qu & (sy\sib y), \\
qr^\prec_{sy,sw}(q) + (q-1)r^\prec_{sy,w}(q) \qu & (y \sib sy).
\end{cases} \nonumber
\end{align}
\item[$(2)$] It holds that \begin{align}
t^\prec_{y,w}(q) &= t^\prec_{sy,sw}(q). \nonumber
\end{align}
\end{enumerate}
\end{prop}

\begin{proof}
$(1)$ If $sy \sib y$, then by Corollary \ref{10.3.3}, there exists a bijection $\psi:P^\prec_r(sy,sw) \rightarrow P^\prec_r(y,w)$ such that $\ell(\psi(\nb)) = \ell(\nb)$ for all $\nb \in P^\prec_r(sy,sw)$. Also, we have $\ell^\si(y,w) = \ell^\si(sy,sw)$. Therefore, we see that
\begin{align}
r^\prec_{sy,sw}(q) &= \sum_{\nb \in P^\prec_r(sy,sw)} q^{\frac{1}{2}\bigl(\ell^\si(sy,sw)-\ell(\nb)\bigr)}(q-1)^{\ell(\nb)} \nonumber\\
&= \sum_{\nb' \in P^\prec_r(y,w)} q^{\frac{1}{2}\bigl(\ell^\si(y,w)-\ell(\nb')\bigr)}(q-1)^{\ell(\nb')} \nonumber\\
&= r^\prec_{y,w}(q). \nonumber
\end{align}
Next, assume that $y \sib sy$. By Corollary \ref{10.3.3}, there exists a bijection $\varphi:P^\prec_r(sy,sw) \sqcup P^\prec_r(sy,w) \rightarrow P^\prec_r(y,w)$ such that $\ell(\varphi(\nb)) = \ell(\nb)$ if $\nb \in P^\prec_r(sy,sw)$, and such that $\ell(\varphi(\nb)) = \ell(\nb)+1$ if $\nb \in P^\prec_r(sy,w)$. Also, we have $\ell^\si(y,w) = \ell^\si(sy,w) + 1 = \ell^\si(sy,sw) + 2$. Therefore, we see that
\begin{align}
&qr^\prec_{sy,sw}(q) + (q-1)r^\prec_{sy,w}(q) \nonumber\\
&= \sum_{\nb \in P^\prec_r(sy,sw)} q^{\frac{1}{2}\bigl(\ell^\si(sy,sw)-\ell(\nb)\bigr)+1}(q-1)^{\ell(\nb)} + \sum_{\nb \in P^\prec_r(sy,w)} q^{\frac{1}{2}\bigl(\ell^\si(sy,s)-\ell(\nb)\bigr)}(q-1)^{\ell(\nb)+1} \nonumber\\
&= \sum_{\substack{\nb' \in P^\prec_r(y,w)\\ \varphi\inv(\nb') \in P^\prec_r(sy,sw)}} q^{\frac{1}{2}\bigl(\ell^\si(y,w)-\ell(\nb')\bigr)}(q-1)^{\ell(\nb')} + \sum_{\substack{\nb' \in P^\prec_r(y,w)\\ \varphi\inv(\nb') \in P^\prec_r(sy,w)}} q^{\frac{1}{2}\bigl(\ell^\si(y,w)-\ell(\nb')\bigr)}(q-1)^{\ell(\nb')} \nonumber\\
&= \sum_{\nb' \in P^\prec_r(y,w)} q^{\frac{1}{2}\bigl(\ell^\si(y,w)-\ell(\nb')\bigr)}(q-1)^{\ell(\nb')} \nonumber\\
&= r^\prec_{y,w}(q), \nonumber
\end{align}
as desired.

$(2)$ Let $\nb = (y = u_0 \xrightarrow[\beta_1]{m_1} u_1 \xrightarrow[\beta_2]{m_2} \cdots \xrightarrow[\beta_k]{m_k} u_k = w) \in P^\prec_t(y,w)$. Since every edge in $\nb$ is a translation edge, we have $\ell^\si(u_i,u_{i+1}) = \ell^\si(su_i,su_{i+1})$ for all $i = 0,1,\ldots,k-1$. Hence we obtain a path $(sy = su_0 \xrightarrow[\beta_1]{m_1} su_1 \xrightarrow[\beta_2]{m_2} \cdots \xrightarrow[\beta_k]{m_k} su_k = sw) \in P^\prec_t(sy,sw)$. This procedure obviously gives a bijection from $P^\prec_t(y,w)$ to $P^\prec_t(sy,sw)$, which preserves both of degree and length of paths. Therefore, we see that
\begin{align}
t^\prec_{y,w}(q) &= \sum_{\nb \in P^\prec_t(y,w)} q^{\deg(\nb)}(q-1)^{\ell(\nb)} \nonumber\\
&= \sum_{\nb' \in P^\prec_t(sy,sw)} q^{\deg(\nb')}(q-1)^{\ell(\nb')} \nonumber\\
&= t^\prec_{sy,sw}(q), \nonumber
\end{align}
as desired.
\end{proof}

\begin{defi}\normalfont
For $y,w \in \Waf$, we define
\begin{align}
\mathcal{R}^\prec_{y,w}(q) := \sum_{\nb \in P^\prec(y,w)} q^{\deg(\nb)}(q-1)^{\ell(\nb)}. \nonumber
\end{align}
\end{defi}

In order to prove Theorem \ref{main}, it suffices to show that the $\rp^\prec_{y,w}(q)$'s satisfy the following equations.
\begin{enumerate}
\item[(r1)] $\rp^\prec_{w,w}(q) = 1$ for all $w \in \Waf$.
\item[(r2)] $\rp^\prec_{y,w}(q) = 0$ unless $y \leq_\si w$ for all $y,w \in \Waf$.
\item[(r3)] $\rp^\prec_{y,w}(q) = \begin{cases}
\rp^\prec_{sy,sw}(q) \qu & \mathrm{if}\ sw\sib w\ \mathrm{and}\ sy\sib y, \\
q\rp^\prec_{sy,sw}(q) +(q-1)\rp^\prec_{y,sw}(q) \qu & \mathrm{if}\ sw\sib w\ \mathrm{and}\ y\sib sy, \end{cases}$\\
for all $y,w \in \Waf$ and all $s \in S_{\af}$.
\item[(r4)] $\sum_{x \in W} (-1)^{\ell^\si(y,xt_\lambda)}q^{\ell(xw_0)}\rp^\prec_{y,xt_\lambda}(q) = \delta_{\lambda,\lambda(y)}$ for all $\lambda \in \Qv$, $y \in \Waf$.
\end{enumerate}
Equations (r1) and (r2) are obvious from the definition of $\rp^\prec_{y,w}(q)$'s. In the next two subsections, we prove the remaining equations.

\subsection{Proof of equation (r3)}
In order to prove equation (r3), we need the equation
\begin{align}
\rp^\prec_{y,w}(q) = \sum_{\mu \in \Qv_+} r^<_{y,wt_{-\mu}}(q)t^\prec_{wt_{-\mu},w}(q); \label{eq:8}
\end{align}
by using this, equation (r3) follows immediately from Proposition \ref{11.1.2} (see Corollary \ref{11.2.3}). Note that the right-hand side of equation \eqref{eq:8} is actually a finite sum because by the definition of $r^\prec_{y,wt_{-\mu}}(q)$, it is equal to zero unless $y \sib wt_{-\mu}$, and because there are only finitely many $\mu \in \Qv_+$ satisfying $y \sib wt_{-\mu}$. The idea for the proof of this equation is to divide a path $\nb \in P^\prec(y,w)$ into two paths $\nb_1 \in P^\prec_r(y,wt_{-\mu})$ and $\nb_2 \in P^\prec_t(wt_{-\mu},w)$ for each $\mu \in \Qv_+$ as follows.

Let us take an arbitrary path $\nb = (y = u_0 \xrightarrow[\beta_1]{m_1} \cdots \xrightarrow[\beta_k]{m_k} u_{k} = w) \in P^\prec(y,w)$, and write $T(\nb) = (m_\beta)_{\beta \in \Phi_+}$. We define $(m'_\beta)_{\beta \in \Phi_+} \in (\Z_{\geq 0})^{|\Phi_+|}$ by:
\begin{align}
m'_\beta = \begin{cases}
m_\beta \qu & \mathrm{if}\ \beta \notin r(\nb), \\
m_i \qu & \mathrm{if}\ \beta = \beta_i \in r(\nb),\ u_i = \cl(u_{i-1})s_{\beta_i} t_{\wt(u_{i-1}) + m_i\beta_i^\vee}\ \mathrm{and}\ \cl(u_{i-1}) < \cl(u_{i-1})s_{\beta_i},\\
m_i -1 \qu & \mathrm{if}\ \beta = \beta_i \in r(\nb),\ u_i = \cl(u_{i-1})s_{\beta_i} t_{\wt(u_{i-1}) + m_i\beta_i^\vee}\ \mathrm{and}\ \cl(u_{i-1})s_{\beta_i} < \cl(u_{i-1}),
\end{cases} \nonumber
\end{align}
and set $\mu' = \sum_{\beta \in \Phi_+} m'_\beta \beta^\vee \in \Qv_+$. Let $(n_\beta)_{\beta \in \Phi_+} \in (\Z_{\geq 0})^{|\Phi_+|}$ be such that
\begin{align}
&n_\beta = m'_\beta \qu (\beta \notin r(\nb)),\ \mathrm{and} \label{assump1}\\ 
&0 \leq n_\beta \leq m'_\beta \qu (\beta \in r(\nb)).\label{assump2}
\end{align}
Then, for $\mu := \sum_{\beta \in \Phi_+} n_\beta \beta^\vee$, there exists a unique pair $(\nb_{(n_\beta),r},\nb_{(n_\beta),t}) \in P^\prec_r(y,wt_{-\mu}) \times P^\prec_t(wt_{-\mu},w)$ such that $r(\nb_{(n_\beta),r}) = r(\nb)$, $T(\nb_{(n_\beta),t}) = (n_\beta)_{\beta \in \Phi_+}$, and $T(\nb) = T(\nb_{(n_\beta),r}) + T(\nb_{(n_\beta),t})$, where the sum is taken componentwise. More precisely, we define $\nb_{(n_\beta),r}$ and $\nb_{(n_\beta),t}$ as follows. If we write $\{ 1 \leq i \leq \ell(\nb) \mid \beta_i \in r(\nb) \} = \{ i_1,i_2,\ldots,i_l \}$ with $i_1 < \cdots < i_l$, and $\{ 1 \leq j \leq \ell(\nb) \mid n_j \neq 0 \} = \{ j_1,j_2,\ldots,j_{l'} \}$ with $j_1 < j_2 < \cdots < j_{l'}$, then
\begin{align}
\nb_{(n_\beta),r} &:= (y \xrightarrow[\beta_{i_1}]{m_{\beta_{i_1}}-n_{\beta_{i_1}}} \cdots \xrightarrow[\beta_{i_l}]{m_{\beta_{i_l}} - n_{\beta_{i_l}}} wt_{-\mu}), \nonumber\\
\nb_{(n_\beta),t} &:= (wt_{-\mu} \xrightarrow[\beta_{j_1}]{n_{\beta_{j_1}}} \cdots \xrightarrow[\beta_{j_{l'}}]{n_{\beta_{j_{l'}}}} w). \nonumber
\end{align}

Conversely, let $(\nb_1,\nb_2) \in P^\prec_r(y,wt_{-\mu}) \times P^\prec_t(wt_{-\mu},w)$ be such that $r(\nb_1) = r(\nb)$ and $T(\nb_1) + T(\nb_2) = T(\nb)$ for some $\mu \in \Qv_+$. We write $T(\nb_1) = (n_{1,\beta})_{\beta \in \Phi_+}$ and $T(\nb_2) = (n_{2,\beta})_{\beta \in \Phi_+}$. Since $\nb_2 \in P^\prec_t(wt_{-\mu},w)$, we have $\mu = \sum_{\beta \in \Phi_+} n_{2,\beta} \beta^\vee$. Also, $T(\nb_1) + T(\nb_2) = T(\nb)$ implies that $(n_{2,\beta})_{\beta \in \Phi_+}$ satisfies condition \eqref{assump1} and that $0 \leq n_\beta \leq m_\beta$ for all $\beta \in r(\nb)$. Moreover, if $\beta = \beta_i \in r(\nb)$, $u_i = \cl(u_{i-1})s_\beta t_{\wt(u_{i-1}) + m\beta^\vee}$, and $\cl(u_{i-1})s_\beta < \cl(u_{i-1})$, then by the definition, we have $n_{1,\beta_i} \geq 1$. Hence $(n_{2,\beta})_{\beta \in \Phi_+}$ satisfies condition \eqref{assump2}. Thefore, $\nb_{(n_{2,\beta}),r}$ and $\nb_{(n_{2,\beta}),t}$ defined above agree with $\nb_1$ and $\nb_2$, respectively.

Thus we obtain a one-to-one correspondence between the set
\begin{align}
A(\nb) := \{ (n_\beta)_{\beta \in \Phi_+} \in (\Z_{\geq 0})^{|\Phi_+|} \mid (n_\beta)_{\beta \in \Phi_+} \textrm{ satisfies conditions } \eqref{assump1}\textrm{ and } \eqref{assump2} \} \nonumber
\end{align}
and the set $B(\nb) := \bigsqcup_{\mu \in \Qv_+} B(\nb)_\mu$, where
\begin{align}
B(\nb)_\mu := \{ (\nb_1,\nb_2) \in P^\prec_r(y,wt_{-\mu}) \times P^\prec_t(wt_{-\mu},w) \mid r(\nb_1) = r(\nb)\ \mathrm{and}\ T(\nb_1) + T(\nb_2) = T(\nb)\}. \nonumber
\end{align}

\begin{lem}\label{11.2.1}
For each $\nb \in P^\prec(y,w)$, we have
\begin{align}
q^{\deg(\nb)}(q-1)^{\ell(\nb)} = \sum_{\mu \in \Qv_+} \sum_{\substack{(\nb_1,\nb_2) \in B(\nb)\\ \nb_1 \in P^\prec_r(y,wt_{-\mu})}} q^{\frac{1}{2}\bigl(\ell^\si(y,wt_{-\mu}) - \ell(\nb_1)\bigr) + \deg(\nb_2)}(q-1)^{\ell(\nb_1) + \ell(\nb_2)}. \nonumber
\end{align}
\end{lem}

\begin{proof}
We employ the notation above. Noting that $(m'_{\beta})_{\beta \in \Phi_+} \in A(\nb)$, we set
\begin{align}
a(q) = a(\nb)(q) := q^{\frac{1}{2}\bigl(\ell^\si(y,wt_{-\mu'}) - \ell(\nb_{(m'_\beta),r})\bigr) + \deg(\nb_{(m'_\beta),t})} (q-1)^{\ell(\nb_{(m'_\beta),r}) + \ell(\nb_{(m'_\beta),t})}. \nonumber
\end{align}
For each $(n_\beta)_{\beta \in \Phi_+} \in A(\nb)$ with $\sum_{\beta \in \Phi_+} n_\beta \beta^\vee = \mu$, we see from the definitions of $\nb_{(n_\beta),r}$ and $\nb_{(n_\beta),t}$ that
\begin{align}
\ell(\nb_{(n_\beta),r}) &= \left|r(\nb) \right|, \label{ob1}\\
\ell(\nb_{(n_\beta),t}) &= \left|\{ \beta \in \Phi_+ \mid n_\beta \neq 0 \}\right|. \label{ob2}
\end{align}
Also, we have
\begin{align}
\deg(\nb_{(n_\beta),t}) &= \sum_{j = 1}^{\ell(\nb_{(n_\beta),t})}d(wt_{-\mu + n_{\beta_{i_1}}\beta_{i_1}^\vee + \cdots + n_{\beta_{i_{j-1}}}\beta_{i_{j-1}}^\vee} \xrightarrow[\beta_{i_j}]{n_{\beta_{i_j}}} wt_{-\mu + n_{\beta_{i_1}}\beta_{i_1}^\vee + \cdots + n_{\beta_{i_{j}}}\beta_{i_{j}}^\vee}) - \ell(\nb_{(n_\beta),t}) \nonumber\\
&= \frac{1}{2}\ell^\si(wt_{-\mu},w) + \sum_{\beta \in \Phi_+} n_\beta - \ell(\nb_{(n_\beta),t}). \label{ob3}
\end{align}
If we write $r(\nb) \cap \{ \beta \in \Phi_+ \mid m'_\beta \neq 0 \} = \{ \beta_{i_1},\ldots,\beta_{i_{k'}} \}$ and $m'_j = m'_{\beta_{i_j}}$ for $1 \leq j \leq k'$, then we compute as follows:
\begin{align}
\sum_{\mu \in \Qv_+} &\sum_{(\nb_1,\nb_2) \in B(\nb)_\mu } q^{\frac{1}{2}\bigl(\ell^\si(y,wt_{-\mu}) - \ell(\nb_1)\bigr) + \deg(\nb_2)}(q-1)^{\ell(\nb_1) + \ell(\nb_2)} \nonumber\\
&= \sum_{\mu \in \Qv_+} \sum_{\substack{(n_\beta) \in A(\nb)\\ \sum_{\beta} n_\beta \beta^\vee = \mu}} q^{\frac{1}{2}\bigl(\ell^\si(y,wt_{-\mu}) - \ell(\nb_{(n_\beta),r})\bigr) + \deg(\nb_{(n_\beta),t})}(q-1)^{\ell(\nb_{(n_\beta),r}) + \ell(\nb_{(n_\beta),t})} \nonumber\\
&= a(q) \sum_{\mu \in \Qv_+} \sum_{\substack{(n_\beta) \in A(\nb)\\ \sum_{\beta} n_\beta \beta^\vee = \mu}} q^{\frac{1}{2}\ell^\si(wt_{-\mu'},wt_{-\mu}) + \deg(\nb_{(n_\beta),t}) - \deg(\nb_{(m'_\beta),t})} (q-1)^{\ell(\nb_{(n_\beta),t}) - \ell(\nb_{(m'_\beta),t})} \nonumber\\
&= a(q) \sum_{\mu \in \Qv_+} \sum_{\substack{(n_\beta) \in A(\nb)\\ \sum_{\beta} n_\beta \beta^\vee = \mu}} q^{-\sum_\beta (m'_\beta - n_\beta) -\ell(\nb_{(n_\beta),t}) + \ell(\nb_{(m'_\beta),t})} (q-1)^{\ell(\nb_{(n_\beta),t}) - \ell(\nb_{(m'_\beta),t})} \qu (\mathrm{by}\ \eqref{ob3})\nonumber\\
&= a(q) \sum_{\mu \in \Qv_+} \sum_{\substack{(n_\beta) \in A(\nb)\\ \sum_{\beta} n_\beta \beta^\vee = \mu}} q^{-\sum_\beta (m'_\beta - n_\beta)} \left(\frac{q}{q-1} \right)^{-\ell(\nb_{(n_\beta),t}) + \ell(\nb_{(m'_\beta),t})} \nonumber\\
&= a(q) \sum_{n_1 = 0}^{m'_{1}}\sum_{n_2 = 0}^{m'2} \cdots \sum_{n_{k'} = 0}^{m'_{k'}} q^{-\sum_{j = 1}^{k'} (m'_{j} - n_j)} \left(\frac{q}{q-1} \right)^{k' - \left|\{ 1 \leq j \leq k'\ |\  n_j \neq 0 \}\right|} \nonumber\\
&= a(q) \sum_{l_1 = 0}^{m'_{1}}\sum_{l_2 = 0}^{m'_{2}} \cdots \sum_{l_{k'} = 0}^{m'_{k'}} q^{-\sum_{j = 1}^{k'} l_j} \left(\frac{q}{q-1} \right)^{k' - \left|\{1 \leq j \leq k'\ |\ l_j \neq m_{j} \}\right|} \qu (\mathrm{by\ changing}\ n_j\ \mathrm{with}\ l_j) \nonumber\\
&= a(q) \sum_{l_1 = 0}^{m'_{1}}\sum_{l_2 = 0}^{m'_{2}} \cdots \sum_{l_{k'} = 0}^{m'_{k'}} q^{-\sum_{j = 1}^{k'} l_j} \left(\frac{q}{q-1} \right)^{\left|\{ 1 \leq j \leq k'\ |\ n_j = m_{j} \}\right|} \nonumber\\
&= a(q)\left(1 + \sum_{j = 1}^{k'} \sum_{1 \leq i_1 < \cdots < i_j \leq k'} \sum_{l_{i_1} = 1}^{m'_{i_1}} \cdots \sum_{l_{i_j} = 1}^{m'_{i_j}} q^{-l_{i_1} - \cdots - l_{i_j}} \left(\frac{q}{q-1} \right)^{\left|\{ 1 \leq r \leq j\ |\ l_{i_r} = m'_{i_r} \} \right|}\right). \nonumber
\end{align}

Here, by induction on $j$, we deduce that
\begin{align}
\sum_{l_{i_1} = 1}^{m'_{i_1}} \cdots \sum_{l_{i_j} = 1}^{m'_{i_j}} q^{-l_{i_1} - \cdots - l_{i_j}} \left(\frac{q}{q-1} \right)^{\left|\{ 1 \leq r \leq j\ |\ l_{i_r} = m'_{i_r}\}\right|} = \left(\frac{1}{q-1} \right)^j. \nonumber
\end{align}
Indeed, when $j=1$, we have
\begin{align}
\sum_{l_{i_1} = 1}^{m'_{i_1}}q^{-l_{i_1}} \left(\frac{q}{q-1} \right)^{\delta_{l_{i_1},m'_{i_1}}} &= q^{-m'_{i_1}}\frac{q}{q-1} + \sum_{l_{i_1} = 1}^{m'_{i_1}-1}q^{-l_{i_1}} \nonumber\\
&= \frac{q^{-m'_{i_1}+1}}{q-1} + \frac{q\inv(1 - q^{-m'_{i_1}+1})}{1 - q\inv} \nonumber\\
&= \frac{q^{-m'_{i_1}+1}}{q-1} + \frac{1-q^{-m'_{i_1}+1}}{q-1} = \frac{1}{q-1}. \nonumber
\end{align}
When $j \geq 2$, we obtain
\begin{align}
\sum_{l_{i_1} = 1}^{m'_{i_1}} \cdots &\sum_{l_{i_j} = 1}^{m'_{i_j}} q^{-l_{i_1} - \cdots - l_{i_j}} \left(\frac{q}{q-1} \right)^{\left|\{ 1 \leq r \leq j\ |\ l_{i_r} = m'_{i_r}\}\right|} \nonumber\\
&= \sum_{l_{i_2} = 1}^{m'_{i_1}} \cdots \sum_{l_{i_j} = 1}^{m'_{i_j}} q^{-l_{i_2} - \cdots - l_{i_j}} \left(\frac{q}{q-1} \right)^{\left|\{ 2 \leq r \leq j\ |\ l_{i_r} = m'_{i_r}\}\right|} \sum_{l_{i_1} = 1}^{m'_{i_1}}q^{-l_{i_1}}\left(\frac{q}{q-1} \right)^{\delta_{l_{i_1},m'_{i_1}}} \nonumber\\
&= \sum_{l_{i_2} = 1}^{m'_{i_1}} \cdots \sum_{l_{i_j} = 1}^{m'_{i_j}} q^{-l_{i_2} - \cdots - l_{i_j}} \left(\frac{q}{q-1} \right)^{\left|\{ 2 \leq r \leq j\ |\ l_{i_r} = m'_{i_r}\}\right|}\frac{1}{q-1} \nonumber\\
&= \left(\frac{1}{q-1} \right)^{j-1}\frac{1}{q-1} \qu (\mathrm{by\ our\ induction\ hypothesis}) \nonumber\\
&= \left(\frac{1}{q-1} \right)^j. \nonumber
\end{align}

Therefore, we deduce that
\begin{align}
a(q)&\left(1 + \sum_{j = 1}^{k'} \sum_{1 \leq i_1 < \cdots < i_j \leq k'} \sum_{l_{i_1} = 1}^{m'_{i_1}} \cdots \sum_{l_{i_j} = 1}^{m'_{i_j}} q^{-l_{i_1} - \cdots - l_{i_j}} \left(\frac{q}{q-1} \right)^{\left|\{ 1 \leq r \leq j\ |\ l_{i_r} = m'_{i_r}\}\right|}\right) \nonumber\\
&= a(q)\left(1 + \sum_{j=1}^{k'}\sum_{1 \leq i_1 < \cdots < i_j \leq k'} \left(\frac{1}{q-1} \right)^j\right) \nonumber\\
&= a(q)\left(\sum_{j = 1}^{k'}\binom{k'}{j}\left(\frac{1}{q-1} \right)^j\right) \nonumber\\
&= a(q)\left(\frac{q}{q-1}\right)^{k'} \nonumber\\
&= q^{\frac{1}{2}\bigl(\ell^\si(y,wt_{-\mu'}) - \ell(\nb_{(m'_\beta),r})\bigr) + \deg(\nb_{(m'_\beta),t}) + k'}(q-1)^{\ell(\nb_{(m'_\beta),r}) + \ell(\nb_{(m'_\beta),t}) - k'} \qu (\mathrm{by\ the\ definition\ of}\ a(q)). \nonumber
\end{align}
To finish the proof of the lemma, we need to show the following:
\begin{align}
\frac{1}{2}\bigl(\ell^\si(y,wt_{-\mu'}) - \ell(\nb_{(m'_\beta),r})\bigr) + \deg(\nb_{(m'_\beta),t}) + k' &= \deg(\nb), \nonumber\\
\ell(\nb_{(m'_\beta),r}) + \ell(\nb_{(m'_\beta),t}) - k' &= \ell(\nb). \nonumber
\end{align}
These are verified as follows:
\begin{align}
\ell(&\nb_{(m'_\beta),r}) + (\ell(\nb_{(m'_\beta),t}) - k') = |r(\nb)| + |t(\nb)| = \ell(\nb) \qu (\mathrm{by}\ \eqref{ob1}\ \mathrm{and}\ \eqref{ob2}), \nonumber\\
\frac{1}{2}(&\ell^\si(y,wt_{-\mu'}) - \ell(\nb_{(m'_\beta),r})) + \deg(\nb_{(m'_\beta),t}) + k' \nonumber\\
&= \frac{1}{2}\ell^\si(y,w) - \frac{1}{2}\ell(\nb_{(m'_\beta),r}) + \sum_{\beta \in \Phi_+}m'_\beta -\bigl(\ell(\nb_{(m'_\beta),t}) - k'\bigr) \nonumber\\
&= \frac{1}{2}\ell^\si(y,w) + \frac{1}{2}\ell(\nb_{(m'_\beta),r}) + \sum_{\beta \in \Phi_+}m'_\beta -\bigl(\ell(\nb_{(m'_\beta),r}) + \ell(\nb_{(m'_\beta),t}) - k'\bigr) \nonumber\\
&= \frac{1}{2}\ell^\si(y,w) + \frac{1}{2}\ell(\nb_{(m'_\beta),r}) + \sum_{\beta \in \Phi_+}m'_\beta -\ell(\nb), \nonumber\\
\deg(\nb) &= \frac{1}{2}\ell^\si(y,w) + \sum_{\beta \in \Phi_+}m_\beta + \frac{1}{2}\left|\{ \beta \in r(\nb) \mid m'_\beta = m_\beta \}\right| - \frac{1}{2}\left|\{ \beta \in r(\nb) \mid m'_\beta = m_\beta - 1\}\right| - \ell(\nb) \nonumber\\
&= \frac{1}{2}\ell^\si(y,w) + \sum_{\beta \in \Phi_+}m_\beta + \frac{1}{2}|r(\nb)| - \left|\{ \beta \in r(\nb) \mid m'_\beta = m_\beta - 1 \}\right| - \ell(\nb) \nonumber\\
&= \frac{1}{2}\ell^\si(y,w) + \sum_{\beta \in \Phi_+}m'_\beta + \frac{1}{2}\ell(\nb_{(m'_\beta),r}) - \ell(\nb). \nonumber
\end{align}
This completes the proof of the lemma.
\end{proof}

\begin{prop}\label{11.2.2}
For each $y,w \in \Waf$, we have
\begin{align}
\rp^\prec_{y,w}(q) = \sum_{\mu \in \Qv_+}r^\prec_{y,wt_{-\mu}}(q)t^\prec_{wt_{-\mu},w}(q). \nonumber
\end{align}
\end{prop}

\begin{proof}
We compute as follows:
\begin{align}
\mathcal{R}^\prec_{y,w}(q) &= \sum_{\nb \in P^\prec(y,w)} q^{\deg(\nb)}(q-1)^{\ell(\nb)} \nonumber\\
&= \sum_{\nb \in P^\prec(y,w)} \sum_{\mu \in \Qv_+} \sum_{(\nb_1,\nb_2) \in B(\nb)_\mu} q^{\frac{1}{2}\bigl(\ell^\si(y,wt_{-\mu}) - \ell(\nb_1)\bigr) + \deg(\nb_2)}(q-1)^{\ell(\nb_1) + \ell(\nb_2)} \qu (\mathrm{by\ Lemma}\ \ref{11.2.1})\nonumber\\
&= \sum_{\mu \in \Qv_+} \sum_{\nb_1 \in P^\prec_r(y,wt_{-\mu})} \sum_{\nb_2 \in P^\prec_t(wt_{-\mu},w)} q^{\frac{1}{2}\bigl(\ell^\si(y,wt_{-\mu}) - \ell(\nb_1)\bigr) + \deg(\nb_2)}(q-1)^{\ell(\nb_1) + \ell(\nb_2)} \nonumber\\
&= \sum_{\mu \in \Qv_+} \sum_{\nb_1 \in P^\prec_r(y,wt_{-\mu})} q^{\frac{1}{2}\bigl(\ell^\si(y,wt_{-\mu}) - \ell(\nb_1)\bigr)} (q-1)^{\ell(\nb_1)} \sum_{\nb_2 \in P^\prec_t(wt_{-\mu},w)} q^{\deg(\nb_2)} (q-1)^{\ell(\nb_2)} \nonumber\\
&= \sum_{\mu \in \Qv_+}r^\prec_{y,wt_{-\mu}}(q)t^\prec_{wt_{-\mu},w}(q) \qu (\mathrm{by\ the\ definitions}). \nonumber
\end{align}
This proves the proposition.
\end{proof}

\begin{cor}\label{11.2.3}
For each $y,w \in \Waf$, we have
\begin{align}
\rp^\prec_{y,w}(q) = \begin{cases}
\rp^\prec_{sy,sw}(q) \qu & \mathrm{if}\ sw\sib w\ \mathrm{and}\ sy \sib y, \\
q\rp^\prec_{sy,sw}(q) +(q-1)\rp^\prec_{y,sw}(q) \qu & \mathrm{if}\ sw \sib w\ \mathrm{and}\ y\sib sy.
\end{cases}\nonumber
\end{align}
\end{cor}
\begin{proof}
By Proposition \ref{11.1.2} (2), we have
\begin{align}
t^\prec_{xt_{-\mu},x}(q) = t^\prec_{wt_{-\mu},w}(q) \nonumber
\end{align}
for all $x,w \in \Waf$ and all $\mu \in \Qv$. Hence the assertion follows from equation \eqref{eq:8} and Proposition \ref{11.1.2}.
\end{proof}

\subsection{Proof of equation (r4)}
In this subsection, we prove equation (r4):
\begin{align}
\sum_{x \in W} (-1)^{\ell^\si(y,xt_\lambda)}q^{\ell(xw_0)}\rp^\prec_{y,xt_\lambda}(q) = \delta_{\lambda,\lambda(y)}\ \ (\lambda \in \Qv,\ y \in \Waf).\nonumber
\end{align}
By translating all elements in a path by $\lambda \in \Qv$, we obtain a bijection $P^\prec(y,w) \rightarrow P^\prec(yt_\lambda,wt_\lambda)$ for all $y,w \in \Waf$. Hence we may assume that $y \in W$ (or equivalently, $\lambda(y) = 0$). In this subsection, we write $\Phi_+ = \{ \beta_1,\ldots,\beta_n \}$, with $\beta_1 \prec \cdots \prec \beta_n$ ; here, $\preceq$ is the fixed reflection order. To prove equation (r4), we will prove stronger statements:
\begin{align}
\sum_{x \in W} (-1)^{\ell^\si(y,x)}q^{\ell(xw_0)}\rp^\prec_{y,x}(q) &= 1, \nonumber\\
\sum_{x \in W} (-1)^{\ell^\si(y,xt_\lambda)} q^{\ell(xw_0)} \sum_{\substack{\nb \in P^\prec(y,xt_\lambda) \\ \beta_k \preceq e(\nb)}} q^{\deg{\nb}}(q-1)^{\ell(\nb)} &= 0\ \ (\lambda \in \Qv\setminus\{0\},\ k \in \{ 1,\ldots,n \}),\nonumber
\end{align}
where $\beta_k \preceq e(\nb)$ means that $\beta_k \preceq \beta$ for all $\beta \in e(\nb)$.

\begin{prop}\label{finite}
If $w \in W$, then $\rp^\prec_{y,w}(q) = \rp_{y,w}(q)$ for all $y \in W$. In particular, equation $(r4)$ holds when $\lambda = 0$.
\end{prop}
\begin{proof}
The polynomials $\rp_{y,w}(q)$ for $y,w \in W$ are uniquely determined by equations (R1), (R2), and (R3) (see Proposition \ref{periodicR}). Here, recall that we have proved that the polynomials $\rp^\prec_{y,w}(q)$ for $y,w \in W$ satisfy the corresponding equations (r1), (r2), and (r3). Therefore, it follows that $\rp^\prec_{y,w}(q) = \rp_{y,w}(q)$ for all $y,w \in W$. In particular, equation (r4) with $\lambda = 0$ is just equation (R4) with $\lambda = 0$. Hence the second assertion follows from Proposition \ref{periodicR}.
\end{proof}

\begin{lem}\label{11.3.2}
For every $z \in W$ and $k \in \{1,\ldots,n\}$, we have
\begin{align}
\sum_{x \in W}(-1)^{\ell(z,x)}q^{\ell(xw_0)}\sum_{\substack{\nb \in P^\prec(z,x) \\ \beta_k \preceq e(\nb)}} q^{\deg(\nb)}(q-1)^{\ell(\nb)} = q^{\left| \{ j \ |\ j<k\ \mathrm{and}\ z<zs_{\beta_j}\} \right|}. \nonumber
\end{align}
\end{lem}
\begin{proof}
We prove the desired equation by descending induction on $\ell(z)$. When $\ell(z) = \ell(w_0)$, then $z = w_0$ and hence the desired equation follows immediately.

Now, we assume that $\ell(z) < \ell(w_0)$, and proceed also by induction on $k$. When $k=1$, we have
\begin{align}
L.H.S. &= \sum_{x \in W} (-1)^{\ell(z,x)}q^{\ell(xw_0)} \sum_{\nb \in P^\prec(z,x)} q^{\deg(\nb)}(q-1)^{\ell(\nb)} \nonumber\\
&= \sum_{x \in W} (-1)^{\ell(z,x)}q^{\ell(xw_0)} \rp^\prec_{z,x}(q) =1 \qu (\mathrm{by\ (r4)\ with}\ \lambda=0), \qu \mathrm{and} \nonumber\\
R.H.S. &= q^0 = 1. \nonumber
\end{align}
Hence we assume that $k \geq 2$. Here observe that since $z,x \in W$, every edge of a path in $P^\prec(z,x)$ is a reflection edge, and hence $m_{\beta_l} = 0$ for all $l = 1,\ldots,n$.

If $zs_{\beta_{k-1}} < z$, then there does not exist any path beginning with an edge of the form $z \xrightarrow[\beta_{k-1}]{0} \ $. Therefore, for each $\nb \in P^\prec(z,x)$, we have $\beta_k \preceq e(\nb)$ if and only if $\beta_{k-1} \preceq e(\nb)$. From this, we deduce that
\begin{align}
L.H.S. &= \sum_{x \in W} (-1)^{\ell(z,x)}q^{\ell(xw_0)} \sum_{\nb \in P^\prec(z,x),\ \beta_{k-1} \preceq e(\nb)} q^{\deg(\nb)}(q-1)^{\ell(\nb)} \nonumber\\
&=q^{\left| \{ j\ |\ j < k-1\ \mathrm{and}\ z < zs_{\beta_j}\} \right|} \qu (\mathrm{by\ our\ induction\ hypothesis\ applied\ to} \ k-1)\nonumber\\
&= q^{\left| \{ j\ |\ j<k\ \mathrm{and}\ z < zs_{\beta_j}\} \right|}. \nonumber
\end{align}

If $z < zs_{\beta_{k-1}}$, then we see that
\begin{align}
L.H.S. &= \sum_{x \in W}(-1)^{\ell(z,x)}q^{\ell(xw_0)}\bigl(\sum_{\substack{\nb \in P^\prec(z,x)\\ \beta_{k-1} \preceq e(\nb)}}q^{\deg(\nb)}(q-1)^{\ell(\nb)} - \hspace{-1cm}\sum_{\substack{\nb \in P^\prec(z,x) \\ \beta_{k-1} \preceq e(\nb)\ \mathrm{and}\ \beta_{k-1} \in e(\nb)}}q^{\deg(\nb)}(q-1)^{\ell(\nb)}\bigr) \nonumber\\
&= q^{\left| \{ j\ |\ j<k-1\ \mathrm{and}\ z < zs_{\beta_j}\} \right|} \nonumber\\
&\hspace{1cm}+ \sum_{x \in W}(-1)^{\ell(zs_{\beta_{k-1}},x)}q^{\ell(xw_0)}\sum_{\substack{\nb \in P^\prec(zs_{\beta_{k-1}},x)\\ \beta_{k} \preceq e(\nb)}}q^{\deg(\nb) + \frac{1}{2}\bigl(\ell(z,zs_{\beta_{k-1}}) - 1\bigr)}(q-1)^{\ell(\nb) +1} \nonumber\\
&= q^{\left| \{ j\ |\ j<k-1\ \mathrm{and}\ z<zs_{\beta_j}\} \right|} \nonumber\\
&\hspace{1cm}+ (q-1)q^{\frac{1}{2}\bigl(\ell(z,zs_{\beta_{k-1}})-1\bigr) + \left| \{ j\ |\ j<k\ \mathrm{and}\ zs_{\beta_{k-1}}<\ zs_{\beta_{k-1}}s_{\beta_j}\} \right|}, \nonumber
\end{align}
where the second equality follows from our induction hypothesis applied to $k-1$, and the third one follows from the desired equality for $zs_{\beta_{k-1}}$. Therefore, if we show the equality
\begin{align}
\left| \{ j \mid j<k-1\ \mathrm{and}\ z<zs_{\beta_j}\} \right| = \left| \{ j \mid j<k\ \mathrm{and}\ zs_{\beta_{k-1}} < zs_{\beta_{k-1}}s_{\beta_j}\} \right| + \frac{1}{2}\bigl(\ell(z,zs_{\beta_{k-1}})-1\bigr), \nonumber
\end{align}
then we can complete the proof of the lemma. This equality will be shown in the next lemma.
\end{proof}

\begin{lem}\label{2}
Let $z \in W$, and $k \in \{ 2,\ldots,n \}$ be such that $z < zs_{\beta_{k-1}}$. Then,
\begin{align}
\left| \{ j \mid j<k-1\ \mathrm{and}\ z<zs_{\beta_j}\} \right| = \left| \{ j \mid j<k\ \mathrm{and}\ zs_{\beta_{k-1}} < zs_{\beta_{k-1}}s_{\beta_j}\} \right| + \frac{1}{2}\bigl(\ell(z,zs_{\beta_{k-1}})-1\bigr). \nonumber
\end{align}
\end{lem}

\begin{proof}
Let $i \in \{ 1,\ldots,n \}$. Then we easily deduce the following by using basic facts about the ordinary Bruhat order and length function:
\begin{enumerate}
\item[] $z < zs_{\beta_i}$ if and only if $z(\beta_i) \in \Phi_+$;
\item[] $zs_{\beta_{k-1}} < zs_{\beta_{k-1}}s_{\beta_i}$ if and only if $zs_{\beta_{k-1}}(\beta_i) \in \Phi_+$;
\item[] $\ell(z) = \left|\{j \mid z(\beta_j) \in -\Phi_+ \}\right|$;
\item[] $\ell(zs_{\beta_{k-1}}) = \left|\{ j \mid zs_{\beta_{k-1}}(\beta_j) \in -\Phi_+ \}\right|$.
\end{enumerate}

Now, let $i' \in \{ 1,\ldots,n \}$ be such that $s_{\beta_{k-1}}(\beta_i) \in \{ \beta_{i'},-\beta_{i'} \}$. Then, $i \in \{ 1,\ldots,n \}$ satisfies exactly one of the following:
\begin{enumerate}
\item[{[1]}] $i < k-1$, $z(\beta_i) \in -\Phi_+$, and $zs_{\beta_{k-1}}(\beta_i) \in -\Phi_+$;
\item[{[2]}] $i > k-1$, $z(\beta_i) \in -\Phi_+$, and $zs_{\beta_{k-1}}(\beta_i) \in -\Phi_+$;
\item[{[3]}] $i < k-1$, $z(\beta_i) \in \Phi_+$, and $zs_{\beta_{k-1}}(\beta_i) \in \Phi_+$;
\item[{[4]}] $i > k-1$, $z(\beta_i) \in \Phi_+$, and $zs_{\beta_{k-1}}(\beta_i) \in \Phi_+$;
\item[{[5]}] $i = k-1$;
\item[{[6]}] $i < k-1$, $z(\beta_i) \in -\Phi_+$, $zs_{\beta_{k-1}}(\beta_i) \in \Phi_+$, and $s_{\beta_{k-1}}(\beta_i) \in \Phi_+$;
\item[{[7]}] $i > k-1$, $z(\beta_i) \in -\Phi_+$, $zs_{\beta_{k-1}}(\beta_i) \in \Phi_+$, and $s_{\beta_{k-1}}(\beta_i) \in \Phi_+$;
\item[{[8]}] $i < k-1$, $z(\beta_i) \in \Phi_+$, $zs_{\beta_{k-1}}(\beta_i) \in -\Phi_+$, and $s_{\beta_{k-1}}(\beta_i) \in -\Phi_+$;
\item[{[9]}] $i > k-1$, $z(\beta_i) \in \Phi_+$, $zs_{\beta_{k-1}}(\beta_i) \in -\Phi_+$, and $s_{\beta_{k-1}}(\beta_i) \in -\Phi_+$;
\item[{[10]}] $i < k-1$, $z(\beta_i) \in -\Phi_+$, $zs_{\beta_{k-1}}(\beta_i) \in \Phi_+$, and $s_{\beta_{k-1}}(\beta_i) \in -\Phi_+$;
\item[{[11]}] $i > k-1$, $z(\beta_i) \in -\Phi_+$, $zs_{\beta_{k-1}}(\beta_i) \in \Phi_+$, and $s_{\beta_{k-1}}(\beta_i) \in -\Phi_+$;
\item[{[12]}] $i < k-1$, $z(\beta_i) \in \Phi_+$, $zs_{\beta_{k-1}}(\beta_i) \in -\Phi_+$, and $s_{\beta_{k-1}}(\beta_i) \in \Phi_+$;
\item[{[13]}] $i > k-1$, $z(\beta_i) \in \Phi_+$, $zs_{\beta_{k-1}}(\beta_i) \in -\Phi_+$, and $s_{\beta_{k-1}}(\beta_i) \in \Phi_+$.
\end{enumerate}
In the following, for each $a = 1,2,\ldots,13$, we denote by $c_a$ the number of $l$, $1 \leq l \leq n$, satisfying condition $[a]$. If we are in case $[6]$, then we have $z(s_{\beta_{k-1}}(\beta_i)) \in \Phi_+$, $zs_{\beta_{k-1}}(s_{\beta_{k-1}}(\beta_i)) \in -\Phi_+$, $s_{\beta_{k-1}}(s_{\beta_{k-1}}(\beta_i)) \in \Phi_+$, $s_{\beta_{k-1}}(\beta_i) = \beta_{i'}$, and $s_{\beta_{k-1}}(\beta_i) = \beta_i - \la \beta_i,\beta_{k-1}^\vee \ra \beta_{k-1}$. Since $i < k-1$ in this case, either $i < i' < k-1$ or $i' < i < k-1$ holds. Hence case $[12]$ holds for $i'$. Conversely, by a similar argument, case $[12]$ for $i$ implies case $[6]$ for $i'$. Therefore, we find that $c_6 = c_{12}$. Similarly, we can show that $c_7 = c_{13}$, $c_8 = c_9$, and $c_{10} = c_{11}$.

Finally, from the equalities above, we deduce by using the assertions at the beginning of the proof that
\begin{align}
\left|\{ j < k-1 \mid z < zs_{\beta_j}\}\right| &= c_3 + c_8 + c_{12}, \nonumber\\
\left|\{ j<k \mid zs_{\beta_{k-1}} < zs_{\beta_{k-1}}s_{\beta_j} \}\right| &= c_3 + c_6 + c_{10}, \nonumber\\
\ell(z,zs_{\beta_{k-1}}) &= c_1 + c_2 + c_5 + c_8 + c_9 + c_{12} + c_{13} - (c_1 + c_2 + c_6 + c_7 + c_{10} + c_{11}) \nonumber\\
&= c_5 + 2(c_8-c_{10}), \nonumber\\
c_5 &= 1. \nonumber
\end{align}
These equalities, together with the equality $c_6 = c_{12}$, prove the lemma.
\end{proof}

\begin{prop}\label{3}
For every $z \in W$ and $k \in \{ 1,2,\ldots,n\}$ such that $z < zs_{\beta_k}$, we have
\begin{align}
\sum_{x \in W}(-1)^{\ell(z,x)}q^{\ell(xw_0)}\sum_{\substack{\nb \in P^\prec(z,x)\\ \beta_k \prec e(\nb)}}&q^{\deg(\nb)}(q-1)^{\ell(\nb)} \nonumber\\
 &\hspace{-1cm}+ \frac{q}{q-1}\sum_{x \in W}(-1)^{\ell(z,x)}q^{\ell(xw_0)}\sum_{\substack{\nb \in P^\prec(z,x)\\ \beta_k \preceq e(\nb)\ \mathrm{and}\ \beta_k \in e(\nb)}}q^{\deg(\nb)}(q-1)^{\ell(\nb)}=0. \nonumber
\end{align}
\end{prop}
\begin{proof}
We compute as follows:
\begin{align}
L.H.S. &= \sum_{x \in W}(-1)^{\ell(z,x)} q^{\ell(xw_0)} \sum_{\substack{\nb \in P^\prec(z,x) \\ \beta_k \preceq e(\nb)}}q^{\deg(\nb)} (q-1)^{\ell(\nb)} \nonumber\\
&\hspace{1cm}+ \frac{1}{q-1} \sum_{x \in W} (-1)^{\ell(z,x)}q^{\ell(xw_0)}\sum_{\substack{\nb \in P^\prec(z,x)\\ \beta_k \preceq e(\nb)\ \mathrm{and}\ \beta_k \in e(\nb)}}q^{\deg(\nb)}(q-1)^{\ell(\nb)} \nonumber\\
&= q^{\left|\{ j \mid j<k\ \mathrm{and}\ z < zs_{\beta_j} \}\right|} \nonumber\\
&\hspace{1cm}+ \frac{1}{q-1} \sum_{x \in W} (-1)^{\ell(z,x)}q^{\ell(xw_0)}\sum_{\substack{\nb \in P^\prec(z,x)\\ \beta_k \preceq e(\nb)\ \mathrm{and}\ \beta_k \in e(\nb)}}q^{\deg(\nb)}(q-1)^{\ell(\nb)} \nonumber\\
&= q^{\left|\{ j \mid j<k\ \mathrm{and}\ z < zs_{\beta_j} \}\right|} \nonumber\\
&\hspace{1cm} - \frac{1}{q-1}\sum_{x \in W} (-1)^{\ell(zs_{\beta_k},x)}q^{\ell(xw_0)}\sum_{\substack{\nb' \in P^\prec(zs_{\beta_k},x)\\ \beta_{k+1} \preceq e(\nb')}}q^{\deg(\nb') + \frac{1}{2}\bigl(\ell(z,zs_{\beta_k})-1\bigr)}(q-1)^{\ell(\nb') + 1} \nonumber\\
&= q^{\left|\{ j \mid j<k\ \mathrm{and}\ z < zs_{\beta_j} \}\right|} - q^{\left|\{ j \mid j<k + 1\ \mathrm{and}\ zs_{\beta_k} < zs_{\beta_k}s_{\beta_j} \}\right| + \frac{1}{2}\bigl(\ell(z,zs_{\beta_k})-1\bigr)} \nonumber\\
&= 0 \qu (\mathrm{by\ Lemma}\ \ref{2}), \nonumber
\end{align}
where the second and fourth equalities follow from Lemma \ref{11.3.2}. This proves the proposition.
\end{proof}

\begin{prop}
Let $z \in W$, $\lambda \in \Qv \setminus\{ 0 \}$, and $k \in \{ 1,\ldots,n \}$. Then we have
\begin{align}
\sum_{x \in W} (-1)^{\ell^\si(z,xt_\lambda)}q^{\ell(xw_0)}\sum_{\substack{\nb \in P^\prec(z,xt_\lambda)\\ \beta_k \preceq e(\nb)}} q^{\deg{\nb}}(q-1)^{\ell(\nb)} = 0.\nonumber
\end{align}
\end{prop}

\begin{proof}
If $\lambda \notin \Qv_+$, then the set $P^\prec(z,xt_\lambda)$ is empty since $z \nleq_{\si} xt_\lambda$ for any $x \in W$. Hence we assume that $\lambda \in \Qv_+ \setminus \{0\}$.

Since $\lambda \neq 0$, for each $\nb = (z = z_0 \xrightarrow[\beta_{i_1}]{m_1} \cdots \xrightarrow[\beta_{i_r}]{m_r}z_r = xt_\lambda) \in P^\prec(z,xt_\lambda)$, we can take the largest $j$ for which $z_j \in W$; clearly, $j<r$. With this notation, we define $b(\nb) \in Wt_{\beta_{i_{j+1}}^\vee}$ by:
\begin{align}
b(\nb) = \begin{cases}
z_jt_{\beta_{i_{j+1}}^\vee} \qu & \mathrm{if}\ z_j < z_js_{\beta_{i_{j+1}}},\\
z_js_{\beta_{i_{j+1}}}t_{\beta_{i_{j+1}}^\vee} \qu & \mathrm{if}\ z_j > z_js_{\beta_{i_{j+1}}}.
\end{cases} \nonumber
\end{align}
If we set
\begin{align}
\nb' &= (z = z_0 \xrightarrow[\beta_{i_1}]{m_1} \cdots \xrightarrow[\beta_{i_j}]{m_j} z_j \xrightarrow[\beta_{i_{j+1}}]{1} b(\nb)) \in P^\prec(z,b(\nb)), \nonumber\\
\nb'' &= (b(\nb) = z_{j+1} \xrightarrow[\beta_{i_{j+2}}]{m_{j+2}} \cdots \xrightarrow[\beta_{i_r}]{m_r} z_r = xt_\lambda) \in P^\prec(b(\nb),xt_\lambda), \nonumber\\
\nb''' &= (b(\nb) \xrightarrow[\beta_{i_{j+1}}]{m_{j+1}-1} z_{j+1} \xrightarrow[\beta_{i_{j+2}}]{m_{j+2}} \cdots \xrightarrow[\beta_{i_r}]{m_r} z_r = xt_\lambda) \in P^\prec(b(\nb),xt_\lambda), \nonumber
\end{align}
then we can check by direct calculation that
\begin{align}
q^{\deg(\nb)}(q-1)^{\ell(\nb)} = \begin{cases}
q^{\deg(\nb')}(q-1)^{\ell(\nb')}q^{\deg(\nb'')}(q-1)^{\ell(\nb'')} \qu & \mathrm{if}\ b(\nb) = z_{j+1}, \\
q^{\deg(\nb')}(q-1)^{\ell(\nb')}\frac{q}{q-1}q^{\deg(\nb''')}(q-1)^{\ell(\nb''')} \qu & \mathrm{if}\ b(\nb) \neq z_{j+1}.
\end{cases} \nonumber
\end{align}

Conversely, assume that we are given $b = cl(b)t_{\wt(b)} \in \Waf$, $\nb' \in P^\prec(z,b)$, and $\nb'' \in P^\prec(b,xt_\lambda)$ satisfying the following conditions:
\begin{enumerate}
\item[$(*)$] $\wt(b) = \beta_i^\vee$ for some $i \geq k$;
\item[$(**)$] $\beta_k \preceq e(\nb')$, and all elements of $\Waf$ appearing in $\nb'$ belong to $W$, except for $b$;
\item[$(***)$] $\wt(b) \prec e(\nb'')$.
\end{enumerate}
Then we can construct $\nb \in P^\prec(z,xt_\lambda)$ by concatenating $\nb'$ and $\nb''$; this $\nb$ satisfies
\begin{align}
b(\nb) &= b, \nonumber\\
\beta_k &\preceq e(\nb), \nonumber\\
q^{\deg(\nb)}(q-1)^{\ell(\nb)} &= q^{\deg(\nb')}(q-1)^{\ell(\nb')}q^{\deg(\nb'')}(q-1)^{\ell(\nb'')}. \nonumber
\end{align}
Similarly, if we are given $b \in \Waf$, $\nb' \in P^\prec(z,b)$, and $\nb''' \in P^\prec(b,xt_\lambda)$ satisfying conditions $(*)$, $(**)$, and that $\wt(b) \preceq e(\nb''')$, $\wt(b) \in e(\nb''')$, then we can construct $\nb \in P^\prec(z,xt_\lambda)$ by concatenating $\nb'$ and $\nb'''$; this $\nb$ satisfies
\begin{align}
b(\nb) &= b, \nonumber\\
\beta_k &\preceq e(\nb), \nonumber\\
q^{\deg(\nb)}(q-1)^{\ell(\nb)} &= q^{\deg(\nb')}(q-1)^{\ell(\nb')}\frac{q}{q-1}q^{\deg(\nb''')}(q-1)^{\ell(\nb''')}. \nonumber
\end{align}

From the above, we see that
\begin{align}
&\sum_{x \in W} (-1)^{\ell^\si(z,xt_\lambda)}q^{\ell(xw_0)}\sum_{\substack{\nb \in P^\prec(z,xt_\lambda)\\ \beta_k \preceq e(\nb)}} q^{\deg{\nb}}(q-1)^{\ell(\nb)} \nonumber\\
&= \sum_{\substack{b \in \Waf \\ \beta_k \preceq \wt(b) \in \Phi_+}} \sum_{x \in W}(-1)^{\ell^\si(z,xt_\lambda)} q^{\ell(xw_0)} \sum_{\substack{\nb \in P^\prec(z,xt_\lambda) \\ \beta_k \preceq e(\nb)\ \mathrm{and}\ b(\nb) = b}} q^{\deg(\nb)}(q-1)^{\ell(\nb)} \nonumber\\
&= \sum_{\substack{b \in \Waf \\ \beta_k \preceq \wt(b) \in \Phi_+}} \sum_{\substack{\nb' \in P^\prec(z,b) \\ \nb'\ \mathrm{satisfies}\ (**)}} q^{\deg(\nb')} (q-1)^{\ell(\nb')} \sum_{x \in W} (-1)^{\ell^\si(z,xt_\lambda)} q^{\ell(xw_0)} \nonumber\\
&\cdot \Bigl(\sum_{\substack{\nb'' \in P^\prec(b,xt_\lambda) \\ \wt(b) \prec e(\nb'')}} q^{\deg(\nb'')} (q-1)^{\ell(\nb'')} + \frac{q}{q-1} \sum_{\substack{\nb''' \in P^\prec(b,xt_\lambda) \\ \wt(b) \preceq e(\nb''')\ \mathrm{and}\ \wt(b) \in e(\nb''')}} q^{\deg(\nb''')} (q-1)^{\ell(\nb''')}\Bigr). \nonumber
\end{align}

Now, let us write $\lambda = \sum_{i=1}^l n_i\alpha_i^\vee \in Q^\vee_+$, and set $|\lambda| = \sum_{i=1}^l n_i$. We complete the proof by showing the equation
\begin{align}
\sum_{x \in W} &(-1)^{\ell^\si(z,xt_\lambda)} q^{\ell(xw_0)} \nonumber\\
&\cdot \Bigl(\sum_{\substack{\nb'' \in P^\prec(b,xt_\lambda) \\ \wt(b) \prec e(\nb'')}} q^{\deg(\nb'')} (q-1)^{\ell(\nb'')} + \frac{q}{q-1} \sum_{\substack{\nb''' \in P^\prec(b,xt_\lambda) \\ \wt(b) \preceq e(\nb''')\ \mathrm{and}\ \wt(b) \in e(\nb''')}} q^{\deg(\nb''')} (q-1)^{\ell(\nb''')}\Bigr) = 0 \nonumber
\end{align}
for all $b \in \Waf$ such that $\beta_k \preceq \wt(b) \in \Phi_+$ by induction on $|\lambda|$. When $|\lambda| = 1$, then $\lambda$ is a simple coroot, and hence either $\lambda-\wt(b) \notin \Qv_+$ or $\lambda - \wt(b) = 0$ holds for each $b \in \Waf$ such that $\beta_k \preceq \wt(b) \in \Phi_+$. If $\lambda-\wt(b) \notin \Qv_+$, then the sum in the parentheses above is equal to $0$ since $P^\prec(b,xt_\lambda)$ is empty. If $\lambda - \wt(b) = 0$, then we compute:
\begin{align}
&\sum_{x \in W} (-1)^{\ell^\si(z,xt_\lambda)} q^{\ell(xw_0)}\Bigl(\sum_{\substack{\nb'' \in P^\prec(b,xt_\lambda)\\ \wt(b) \prec e(\nb'')}}q^{\deg(\nb'')}(q-1)^{\ell(\nb'')} + \frac{q}{q-1} \hspace{-1cm} \sum_{\substack{\nb''' \in P^\prec(b,xt_\lambda)\\ \wt(b) \preceq e(\nb''')\ \mathrm{and}\ \wt(b) \in e(\nb''')}} \hspace{-1cm} q^{\deg(\nb''')}(q-1)^{\ell(\nb''')}\Bigr) \nonumber\\
&= \sum_{x \in W} (-1)^{\ell(z,x)} q^{\ell(xw_0)}\Bigl(\sum_{\substack{\nb'' \in P^\prec(\cl(b),x)\\ \wt(b) \prec e(\nb'')}}q^{\deg(\nb'')}(q-1)^{\ell(\nb'')} + \frac{q}{q-1} \hspace{-1cm} \sum_{\substack{\nb''' \in P^\prec(\cl(b),x)\\ \wt(b) \preceq e(\nb''')\ \mathrm{and}\ \wt(b) \in e(\nb''')}} \hspace{-1cm} q^{\deg(\nb''')}(q-1)^{\ell(\nb''')}\Bigr) \nonumber\\
&= (-1)^{\ell(z,\cl(b))}\sum_{x \in W} (-1)^{\ell(\cl(b),x)} q^{\ell(xw_0)} \nonumber\\
&\hspace{2cm}\cdot \Bigl(\sum_{\substack{\nb'' \in P^\prec(\cl(b),x)\\ \wt(b) \prec e(\nb'')}}q^{\deg(\nb'')}(q-1)^{\ell(\nb'')} + \frac{q}{q-1} \hspace{-1cm} \sum_{\substack{\nb''' \in P^\prec(\cl(b),x)\\ \wt(b) \preceq e(\nb''')\ \mathrm{and}\ \wt(b) \in e(\nb''')}} \hspace{-1cm} q^{\deg(\nb''')}(q-1)^{\ell(\nb''')}\Bigr) \nonumber\\
&= 0 \qu (\mathrm{by\ Proposition}\ \ref{3}), \nonumber
\end{align}
where the first equality holds since $(-1)^{\ell^\si(z,xt_\lambda)} = (-1)^{\ell(z,x) + 2\la \rho,\lambda \ra} = (-1)^{\ell(z,x)}$, and the third equality follows from Proposition \ref{3} applied to $\cl(b) \in W$ and $\wt(b) \in \Phi_+$ (note that we have $\cl(b) < \cl(b)s_{\wt(b)}$ by the definitions of $b$ and $b(\nb)$). This proves the desired equation when $|\lambda| = 1$. When $|\lambda| \geq 2$, we compute:
\begin{align}
&\sum_{x \in W} (-1)^{\ell^\si(z,xt_\lambda)} q^{\ell(xw_0)}\Bigl(\sum_{\substack{\nb'' \in P^\prec(b,xt_\lambda)\\ \wt(b) \prec e(\nb'')}}q^{\deg(\nb'')}(q-1)^{\ell(\nb'')} + \frac{q}{q-1} \hspace{-1cm} \sum_{\substack{\nb''' \in P^\prec(b,xt_\lambda)\\ \wt(b) \preceq e(\nb''')\ \mathrm{and}\ \wt(b) \in e(\nb''')}} \hspace{-1cm} q^{\deg(\nb''')}(q-1)^{\ell(\nb''')}\Bigr) \nonumber\\
&= \sum_{x \in W} (-1)^{\ell^\si(z,xt_\lambda)} q^{\ell(xw_0)} \nonumber\\
&\hspace{2cm}\cdot \Bigl( \hspace{-0.25cm} \sum_{\substack{\nb'' \in P^\prec(\cl(b),xt_{\lambda - \wt(b)})\\ \wt(b) \prec e(\nb'')}} \hspace{-0.25cm} q^{\deg(\nb'')}(q-1)^{\ell(\nb'')} + \frac{q}{q-1} \hspace{-1cm} \sum_{\substack{\nb''' \in P^\prec(\cl(b),xt_{\lambda - \wt(b)})\\ \wt(b) \preceq e(\nb''')\ \mathrm{and}\ \wt(b) \in e(\nb''')}} \hspace{-1cm} q^{\deg(\nb''')}(q-1)^{\ell(\nb''')}\Bigr) \nonumber\\
&= 0 \qu (\mathrm{by\ our\ induction\ hypothesis\ applied\ to}\ |\lambda|-\wt(b)). \nonumber
\end{align}
This completes the proof of the proposition.
\end{proof}

\setcounter{theo}{0}
\section{Another description of Theorem \ref{main}}
In the actual computation of periodic $R$-polynomials for which Theorem \ref{main} is used, it is sometimes difficult to determine whether or not a given path is an element of $P^\prec(y,w)$. Hence it would be convenient to describe the paths in $P^\prec(y,w)$ in a different way. For this purpose, we introduce a new finite directed graph, which we call the $double\ Bruhat\ graph\ (\DBG)$. The idea is to take the ``finite" part of a path.

\subsection{The double Bruhat graph}
\begin{defi}\normalfont
The double Bruhat graph associated to a finite Weyl group $W$ and a finite root system $\Phi$ is the directed graph with vertex set $W$, in which two vertices $y'$ and $w'$ are joined by a labeled arrow $y' \xrightarrow[\beta]{d} w'\ (\beta \in \Phi_+,\ d \in \Z_{> 0})$ if and only if one of the following conditions holds:
\begin{enumerate}
\item[$(1)$] $y' < w',\ w' = y's_\beta$ and $d = \frac{1}{2}(\ell(y',w') +1)$;
\item[$(2)$] $y' > w',\ w' = y's_\beta$ and $d = \frac{1}{2}(\ell(y',w') +2\la \rho,\beta^\vee \ra + 1)$.
\end{enumerate}
An edge that satisfy condition $(1)$ $($resp., $(2)$$)$ is called a Bruhat edge $($resp., quantum edge$)$.
\end{defi}

\begin{ex}\normalfont
In type $A_2$, the $\DBG$ is as follows.
\begin{figure}[H]
\centering
\labellist
\large\hair 2pt
\pinlabel $e$ at 169 60
\pinlabel $s_1$ at 375 60
\pinlabel $s_2$ at 70 232
\pinlabel $s_2s_1$ at 475 232
\pinlabel $s_1s_2$ [r] at 190 410
\pinlabel $w_0$ at 375 410
\small
\pinlabel $\alpha_1$ at 270 80
\pinlabel $\alpha_1$ at 360 255
\pinlabel $\alpha_1$ at 270 425
\pinlabel $\alpha_2$ at 138 150
\pinlabel $\alpha_2$ at 245 330
\pinlabel $\alpha_2$ at 443 328
\pinlabel $\alpha_1+\alpha_2$ at 385 164
\pinlabel $\alpha_1+\alpha_2$ at 189 175
\pinlabel $\alpha_1+\alpha_2$ at 78 338
\pinlabel $2$ at 448 135
\pinlabel $2$ at 305 335
\pinlabel $2$ at 136 303
\endlabellist
\includegraphics[width=8cm]{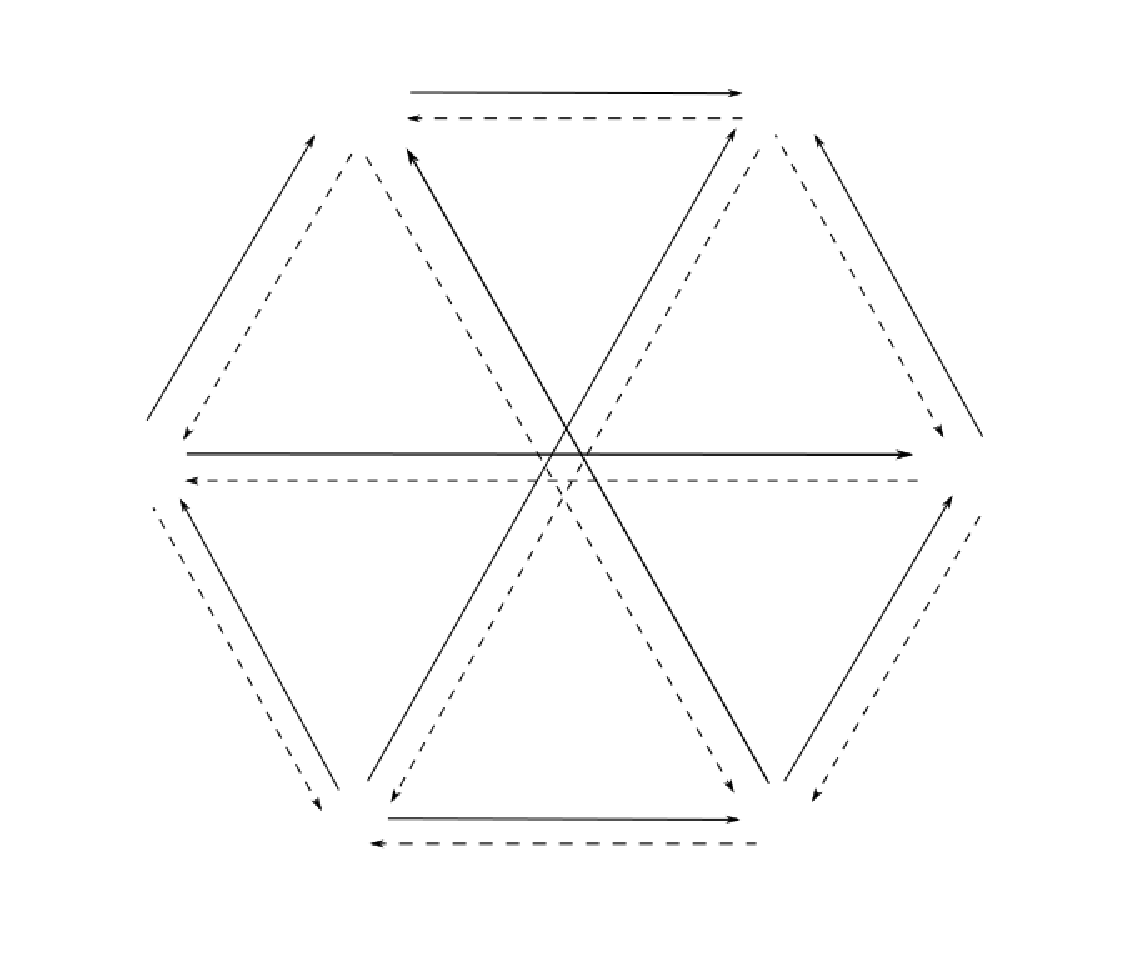}
\end{figure}

\noindent The Bruhat edges and quantum edges are indicated by the ordinary arrows and dotted arrows, respectively; the labels $1$ are omitted.
\end{ex}

Fix a reflection order $\preceq$ on $\Phi_+$. Let $y,w \in \Waf$, and write them as $y = \cl(y)t_{\wt(y)}$ and $w = \cl(w)t_{\wt(w)}$.
\begin{defi}\normalfont
A label-increasing double Bruhat path of length $k$ from $y$ to $w$ with respect to the reflection order $\preceq$ is a sequence $\nb = (\cl(y) = y'_0 \xrightarrow[\beta_1]{d_1} \cdots \xrightarrow[\beta_k]{d_k} y'_k = \cl(w))$ of directed edges in $\DBG$ such that $\beta_1 \preceq \cdots \preceq \beta_k$ and $\sum_{i \in \{1,\ldots,k \ \mid\ y'_{i-1}\xrightarrow[\beta_i]{d_i}y'_i\ \mathrm{is\ quantum} \}} \beta_i^\vee = \wt(w) - \wt(y)$; we denote the length $k$ of a path $\nb$ by $\ell(\nb)$. Let $\DBP^\preceq(y,w)$ denote the set of label-increasing double Bruhat paths from $y$ to $w$.
\end{defi}

\begin{defi}\normalfont
Let $\nb  = (y' = y'_0 \xrightarrow[\beta_1]{d_1} \cdots \xrightarrow[\beta_k]{d_k} y'_k = w') \in \DBP^\preceq(y,w)$.
\begin{enumerate}
\item[$(1)$] We set $e(\nb) := \{ \beta_i \mid i=1,2,\ldots,k\}$; the cardinality of this set is denoted by $\ell'(\nb)$, which does not necessarily agree with $k = \ell(\nb)$.
\item[$(2)$] The degree of $\nb$ is defined to be the integer
\begin{align}
\deg(\nb) := \sum_{i=1}^k d_i -\ell'(\nb). \nonumber
\end{align}
\end{enumerate}
\end{defi}

\subsection{Another description}
\begin{prop}\label{5.0.4}
For each $y,w \in \Waf$, there exists a bijection $\DBP^\preceq(y,w) \rightarrow P^\prec(y,w) ,\nb \mapsto \tilde{\nb}$, such that $\deg(\nb) = \deg(\tilde{\nb})$ and $\ell'(\nb) = \ell(\tilde{\nb})$.
\end{prop}
\begin{proof}
Let $\nb \in \DBP^\preceq(y,w)$, and write it as:
\begin{align}
\nb = (\cl(y) = y'_0 \xrightarrow[\beta_1]{d_{0,1}} y'_{0,1} \xrightarrow[\beta_1]{d_{0,2}} y'_{0,2} \xrightarrow[\beta_1]{d_{0,3}} \cdots \xrightarrow[\beta_1]{d_{0,a_1}} y'_{0,a_1} = y'_1 \xrightarrow[\beta_2]{d_{1,1}} y'_{1,1} \xrightarrow[\beta_2]{d_{1,2}} \cdots \xrightarrow[\beta_k]{d_{k-1,a_k}} y'_{k-1,a_k} = y'_k = \cl(w)), \nonumber
\end{align}
with $\beta_1 \prec \beta_2 \prec \cdots \prec \beta_k$. For each $i =1,\ldots,k$, we define $u_i = \cl(u_i) t_{\wt(u_i)} \in \Waf$ as follows:
\begin{align}
u_0 &= y, \nonumber\\
u_i &= \begin{cases}
u_{i-1} t_{m_i\beta_i^\vee} \qu & (a_i = 2m_i),\\
\cl(u_{i-1})s_{\beta_i}t_{\wt(u_{i-1}) + m_i\beta_i^\vee} \qu & (a_i = 2m_i \pm 1),\\
\end{cases} \nonumber
\end{align}
where $m_i$ is the number of quantum edges in the path $y'_{i-1} \xrightarrow[]{} y'_{i-1,1} \xrightarrow[]{} \cdots \xrightarrow[]{} y'_{i-1,a_i}$; note that if $a_i = 2m_i + 1$, then $y'_{i-1} \xrightarrow[]{} y'_{i-1,1}$ is a Bruhat edge, and if $a_i = 2m_i - 1$ then, $y'_{i-1} \xrightarrow[]{} y'_{i-1,1}$ is a quantum edge. Then, we have $u_{i-1} \xrightarrow[\beta_i]{m_i} u_i$, and
\begin{align}
\deg(u_{i-1} \xrightarrow[\beta_i]{m_i} u_i) = d(u_{i-1} \xrightarrow[\beta_i]{m_i} u_i) - 1 = \begin{cases}
\frac{1}{2}\ell^\si(u_{i-1},u_i) + m_i - 1 \qu & (a_i = 2m_i),\\
\frac{1}{2}\bigl(\ell^\si(u_{i-1},u_i) + 1\bigr) + m_i - 1 \qu & (a_i = 2m_i + 1),\\
\frac{1}{2}\bigl(\ell^\si(u_{i-1},u_i) - 1\bigr) + m_i -1. \qu & (a_i = 2m_i - 1).\\
\end{cases} \nonumber
\end{align}
We now define $\widetilde{\nb}$ to be
\begin{align}
\widetilde{\nb} :=(y = u_0 \xrightarrow[\beta_{1}]{m_1} u_1 \xrightarrow[\beta_2]{m_2} \cdots \xrightarrow[\beta_k]{m_k} u_k = w). \nonumber
\end{align}
It is easy to see that $\widetilde{\nb} \in P^\prec(y,w)$ and $\ell'(\nb) = \ell(\widetilde{\nb})$.

In order to show that $\deg(\nb) = \deg(\widetilde{\nb})$, we need to compare the contribution of $y'_{i-1} \xrightarrow[\beta_i]{d_{i-1,1}} \cdots \xrightarrow[\beta_i]{d_{i-1,a_i}} y'_{i-1,a_i} = y'_i$ to $\deg(\nb)$ and that of $y_{i-1} \xrightarrow[\beta_i]{m_i} y_i$ to $\deg(\widetilde{\nb})$. If $a_i = 2m_i$, then the former is
\begin{align}
\sum_{j=1}^{a_i} d_{i-1,j} - 1, \nonumber
\end{align}
and the latter is
\begin{align}
\deg(u_{i-1} \xrightarrow[\beta_i]{m_i} u_i) &=\frac{1}{2}\ell^\si(u_{i-1},u_i) + m_i -1. \nonumber\\
&=\la\rho,m_i\beta_i^\vee\ra + m_i -1. \nonumber
\end{align}
Here, by the definition of $\DBP^\preceq(y,w)$, we have
\begin{align}
d_{i-1,j-1} + d_{i-1,j} = \la \rho,\beta_i^\vee \ra + 1 \nonumber
\end{align}
for all $j = 2,4,\ldots,a_i$. Therefore, we see that
\begin{align}
\sum_{j=1}^{a_i} d_{i-1,j} - 1 &= \frac{1}{2}a_i(\la \rho,\beta_i^\vee \ra +1) -1 \nonumber\\
&= m_i( \la \rho,\beta_i^\vee \ra +1) -1 \nonumber\\
&= m_i \la \rho,\beta_i^\vee \ra + m_i -1. \nonumber
\end{align}
If $a_i = 2m_i + 1$, then the former is
\begin{align}
d_{i-1,1} + \sum_{j=2}^{a_i} d_{i-1,j} -1, \nonumber
\end{align}
and the latter is
\begin{align}
\deg(u_{i-1} \xrightarrow[\beta_i]{m_i} u_i) &=\frac{1}{2}\bigl(\ell^\si(u_{i-1},u_i) - 1\bigr) + m_i \nonumber\\
&= \frac{1}{2}\bigl(\ell^\si(u_{i-1},\cl(u_{i-1})s_{\beta_i}t_{\wt(u_{i-1})}) + \ell^\si(\cl(u_{i-1})s_{\beta_i}t_{\wt(u_{i-1})},u_i) -1\bigr) + m_i \nonumber\\
&= \frac{1}{2}\bigl(\ell(\cl(u_{i-1}),\cl(u_{i-1})s_{\beta_i}) + 1\bigr) + \la\rho,m_i\beta_i^\vee\ra -1 + m_i \nonumber\\
&= d_{i-1,1} + m_i\la\rho,\beta_i^\vee\ra + m_i -1. \nonumber
\end{align}
In a way similar to the case $a_i = 2m_i$, we deduce that
\begin{align}
\sum_{j=2}^{a_i} d_{i-1,j} = m_i \la \rho,\beta_i^\vee \ra + m_i. \nonumber
\end{align}
If $a_i = 2m_i - 1$, then the former is
\begin{align}
d_{i-1,1} + \sum_{j=2}^{a_i} d_{i-1,j} -1, \nonumber
\end{align}
and the latter is
\begin{align}
\deg(u_{i-1} \xrightarrow[\beta_i]{m_i} u_i) &=\frac{1}{2}\bigl(\ell^\si(u_{i-1},u_i) - 1\bigr) + m_i - 1 \nonumber\\
&= \frac{1}{2}\bigl(\ell^\si(u_{i-1},\cl(u_{i-1})s_{\beta_i}t_{\wt(u_{i-1} + \beta_i^\vee)}) + \ell^\si(\cl(u_{i-1})s_{\beta_i}t_{\wt(u_{i-1}) + \beta_i^\vee},u_i) -1\bigr) + m_i -1 \nonumber\\
&= \frac{1}{2}\bigl(\ell(\cl(u_{i-1}),\cl(u_{i-1})s_{\beta_i}t_{\beta_i^\vee}) + 1\bigr) + \la\rho,(m_i-1)\beta_i^\vee\ra -1 + m_i -1 \nonumber\\
&= d_{i-1,1} + (m_i-1)\la\rho,\beta_i^\vee\ra + m_i -2. \nonumber
\end{align}
In a way similar to the case $a_i = 2m_i$, we deduce that
\begin{align}
\sum_{j=2}^{a_i} d_{i-1,j} = (m_i-1) \la \rho,\beta_i^\vee \ra + m_i-1. \nonumber
\end{align}
Thus, in all cases above, the equality $\deg(\nb) = \deg(\widetilde{\nb})$ holds.

The correspondence $\nb \mapsto \widetilde{\nb}$ is bijective. Indeed, its inverse is given as follows. Let $\widetilde{\nb} = (y = u_0 \xrightarrow[\beta_1]{m_1} \cdots \xrightarrow[\beta_k]{m_k} u_k = w) \in P^\prec(y,w)$. To each $u_{i-1} \xrightarrow[\beta_i]{m_i} u_i$, we assign a path
\begin{align}
(y'_{i-1} \xrightarrow[\beta_i]{d_{i-1,1}} y'_{i-1,1} \xrightarrow[\beta_i]{d_{i-1,2}} \cdots \xrightarrow[\beta_i]{d_{i-1,a_i}} y'_{i-1,a_i} = y'_i), \nonumber
\end{align}
where 
\begin{align}
a_i &= \begin{cases}
2m_i \qu & (u_{i-1} \xrightarrow[\beta_i]{m_i} u_i\ \mathrm{is\ a\ translation\ edge}),\\
2m_i + 1 \qu & (u_{i-1} \xrightarrow[\beta_i]{m_i} u_i\ \mathrm{is\ a\ reflection\ edge\ and}\ \cl(u_{i-1}) < \cl(u_{i-1})s_{\beta_i}),\\
2m_i - 1 \qu & (u_{i-1} \xrightarrow[\beta_i]{m_i} u_i\ \mathrm{is\ a\ reflection\ edge\ and}\ \cl(u_{i-1}) > \cl(u_{i-1})s_{\beta_i}),\\
\end{cases} \nonumber\\
y'_{i-1} &= \cl(u_{i-1}), \nonumber\\
y'_{i-1,j} &= \begin{cases}
\cl(u_{i-1}) \qu & (j\ \mathrm{is\ even}), \\
\cl(u_{i-1})s_{\beta_i} \qu & (j\ \mathrm{is\ odd}).
\end{cases} \nonumber
\end{align}
By concatenating the paths above, we obtain a path $\nb \in \DBP^\preceq(y,w)$ such that $\deg(\nb) = \deg(\widetilde{\nb})$ and $\ell'(\nb) = \ell(\widetilde{\nb})$. This completes the proof of the proposition.
\end{proof}

\begin{theo}
For each $y,w \in \Waf$, we have
\begin{align}
\rp_{y,w}(q) = \sum_{\nb \in \DBP^\preceq(y,w)} q^{\deg(\nb)}(q-1)^{\ell'(\nb)}. \nonumber
\end{align}
\end{theo}
\begin{proof}
This follows from Theorem \ref{main} and Proposition \ref{5.0.4}
\end{proof}

\begin{ex}\normalfont
In type $A_2$, let $y = e$, $w = w_0t_{\alpha_1^\vee+\alpha_2^\vee}$, $\alpha_1 \prec \alpha_1+\alpha_2 \prec \alpha_2$. The elements of $\mathrm{DBP}^\preceq(y,w)$ are as follows:
\begin{align}
\nb_1 &= (e \xrightarrow[\alpha_1+\alpha_2]{2} w_0 \xrightarrow[\alpha_1+\alpha_2]{1} e \xrightarrow[\alpha_1+\alpha_2]{2} w_0), \nonumber\\
\nb_2 &= (e \xrightarrow[\alpha_1]{1} s_1 \xrightarrow[\alpha_1+\alpha_2]{1} s_2s_1 \xrightarrow[\alpha_1+\alpha_2]{2} s_1 \xrightarrow[\alpha_1+\alpha_2]{1} s_2s_1 \xrightarrow[\alpha_2]{1} w_0), \nonumber\\
\nb_3 &= (e \xrightarrow[\alpha_1]{1} s_1 \xrightarrow[\alpha_1]{1} e \xrightarrow[\alpha_1+\alpha_2]{2} w_0 \xrightarrow[\alpha_2]{1} s_2s_1 \xrightarrow[\alpha_2]{1} w_0), \nonumber\\
\nb_4 &= (e \xrightarrow[\alpha_1]{1} s_1 \xrightarrow[\alpha_1]{1} e \xrightarrow[\alpha_1]{1} s_1 \xrightarrow[\alpha_1+\alpha_2]{1} s_2s_1 \xrightarrow[\alpha_2]{1} w_0 \xrightarrow[\alpha_2]{1} s_2s_1 \xrightarrow[\alpha_2]{1} w_0). \nonumber 
\end{align}
In this case, we have
\begin{center}
\begin{tabular}{ll}
$\deg(\nb_1) = 4$, & $\ell'(\nb_1) = 1$, \nonumber\\
$\deg(\nb_2) = 3$, & $\ell'(\nb_2) = 3$, \nonumber\\
$\deg(\nb_3) = 3$, & $\ell'(\nb_3) = 3$, \nonumber\\
$\deg(\nb_4) = 4$, & $\ell'(\nb_4) = 3$. \nonumber
\end{tabular}\end{center}
Therefore, we obtain
\begin{align}
\rp_{e,w_0t_{\alpha_1^\vee+\alpha_2^\vee}}(q) &= q^4(q-1) + 2q^3(q-1)^3 + q^4(q-1)^3. \nonumber
\end{align}
\end{ex}

\end{document}